\documentclass{article}
\usepackage{amsfonts}
\usepackage{graphicx}
\usepackage[latin1]{inputenc}
\usepackage[english]{babel}
\usepackage{color}
\usepackage{float}
\usepackage{fancyhdr}
\usepackage{fancybox}
\usepackage{amssymb,amsmath}
\usepackage{textcomp}
\usepackage{lscape}
\usepackage{amsthm}
\usepackage{parskip}
\usepackage{epstopdf}
\usepackage[labelfont=bf]{caption}
\usepackage{subfigure}
\usepackage[top=90pt, bottom=90pt]{geometry}
\usepackage{longtable}
\usepackage{colortbl}
\usepackage{multirow}
\usepackage[table]{xcolor}
\usepackage{tikz}
\usetikzlibrary{matrix,backgrounds}
\usepackage{algorithm}
\usepackage{algorithmic}
\usepackage{mathabx} 
\numberwithin{equation}{section}
\definecolor{grijs}{RGB}{92,92,92}

\newcommand{\half}{\frac{1}{2}}
\newcommand{\dud}[1]{\frac{\partial v}{\partial #1}}
\newcommand{\dudd}[1]{\frac{\partial^2v}{\partial #1^2}}
\newcommand{\duddm}[2]{\frac{\partial^2v}{\partial #1 \partial #2}}
\newcommand{\dd}{\mathrm{d}}
\newcommand{\rhoh}{\widehat{\rho}}
\newcommand{\oV}{\overline{v}}
\newcommand{\of}{\overline{F}}
\newcommand{\bQ}{\mathbb{Q}}
\newcommand{\bE}{\mathbb{E}}
\newcommand{\minn}{\raisebox{.1pt}{-}}

\definecolor{grijs}{RGB}{232,232,232}

\pgfdeclarelayer{myback}
\pgfsetlayers{myback,background,main}

\tikzset{mycolor/.style = {line width=1bp,color=#1}}
\tikzset{myfillcolor/.style = {draw,fill=#1}}

\author{Lynn Boen\footnote{Department of Mathematics and Computer Science,
University of Antwerp, Middelheimlaan 1, B-2020 Antwerp, Belgium.
\mbox{Email}: \texttt{\{lynn.boen,karel.inthout\}@uantwerpen.be}.}
~and Karel~J.~in 't Hout\footnotemark[\value{footnote}]\\\\
{\it Dedicated to the memory of Willem Hundsdorfer}
}

\title{Operator splitting schemes for the two-asset Merton jump-diffusion model}

\begin{document}
\maketitle
\begin{abstract}
This paper deals with the numerical solution of the two-dimensional time-dependent Merton partial integro-differential 
equation (PIDE) for the values of rainbow options under the two-asset Merton jump-diffusion model.
Key features of this well-known equation are a two-dimensional nonlocal integral part and a mixed spatial derivative 
term.
For its efficient and stable numerical solution, we study seven recent and novel operator splitting schemes of the 
implicit-explicit (IMEX) and the alternating direction implicit (ADI) kind. 
Here the integral part is always conveniently treated in an explicit fashion.
The convergence behaviour and the relative performance of the seven schemes are investigated in ample numerical 
experiments for both European put-on-the-min and put-on-the-average options. 
\end{abstract}


\section{Introduction}
The inherent risk diversification of multi-name or rainbow options, for which the payoff depends on multiple underlying asset prices, 
has made them increasingly popular over the recent years. 
Margrabe~\cite{M78} proposed such a rainbow option, which allows to exchange one asset for another, and derived a valuation formula 
under the Black--Scholes framework. 
The term rainbow option was later introduced by Rubinstein~\cite{R91}, who referred to the number of underlying assets as the number 
of colors of the option. 
The development of fast, accurate and stable methods for obtaining the fair values of options forms a central topic in the field of 
computational finance. 
The fair value of an option is equal to its expected discounted payoff under a risk-neutral measure. 
Acquiring a semi-closed analytical valuation formula, as for example Black \& Scholes~\cite{BS73} and Margrabe~\cite{M78} 
did, would be most favorable, but this is often not possible. 
Accordingly, there is a strong demand for reliable numerical methods that can efficiently approximate fair option values. 

For the numerical approximation of the fair value of an option, three main approaches are considered in the literature up to now.
Firstly, Monte Carlo methods estimate the expected discounted payoff value by computing sample means. 
The second way is by numerical integration, which is employed in for example the Carr--Madan method~\cite{CM99} and the COS 
method by Fang \& Oosterlee~\cite{FO08}. 
The third main approach is to numerically solve a time-dependent partial differential equation (PDE) that holds for the option value.
This approach was introduced in computational finance by Schwartz~\cite{S77}, who considered a finite difference discretization
of the celebrated Black--Scholes PDE. 
The development, analysis and application of numerical methods for time-dependent PDEs has since then been widely studied in the 
literature for valuing many types of options under a broad range of asset pricing models, see for example \cite{H17book,TR00book,W06book}.

It is well-known in financial markets that sudden, large changes in asset prices can occur. 
In view of this common phenomenon, Merton~\cite{M76} added a jump term to the original stochastic process for the asset 
price assumed by Black \& Scholes~\cite{BS73}, yielding a so-called jump-diffusion model.
Under such an advanced asset price process, the option value satisfies a one-dimensional time-dependent \textit{partial 
integro-differential equation} (PIDE). 
The integral part, which stems from the jump term, runs over the whole asset price domain, and hence, is nonlocal. 

In the present paper we consider the extension of the prominent Merton jump-diffusion model to two assets where the jumps 
in the asset prices occur simultaneously.
This yields a two-dimensional time-dependent PIDE for the values of two-asset options. 
The efficient and stable numerical solution of this equation poses two principal challenges: (a) an effective treatment of 
the two-dimensional PDE part, which involves a mixed derivative term, and (b) an effective treatment of the two-dimensional 
integral part, which is nonlocal. 

We follow the well-known and versatile method of lines (MOL) approach \cite{HV03book}, whereby the PIDE is first discretized in the 
spatial variables and subsequently in the temporal variable.
Spatial discretization of the integral part gives rise to a very large, dense matrix that is computationally expensive to work with.
However, by a change of variables, the double integral can be turned into a two-dimensional cross-correlation.  
Spatial discretization of this new integral version leads to a matrix having a block-Toeplitz structure for which matrix-vector 
products can be computed efficiently using the fast Fourier transform (FFT).
For one-dimensional PIDEs, the application of FFT in approximating the integral part has been advocated in \cite{AO05,AA00,HFV05}. 
In particular, Almendral \& Oosterlee \cite{AO05} observed that the matrix corresponding to the jump term is Toeplitz and can be 
embedded in a circulant matrix, so that matrix-vector products can be computed exactly by means of FFT. 
In \cite{CF08,STS14} the FFT technique has next been applied in the case of two-dimensional PIDEs with a two-dimensional 
integral part.
Salmi, Toivanen \& von Sydow \cite{STS14} considered a useful embedding of a block-Toeplitz matrix into a block-circulant 
one.
Here we shall employ a similar idea, which was introduced in a different context in Barrowes, Teixeira \& Kong \cite{BTK01}.

After spatial discretization of the two-dimensional Merton PIDE one arrives at a very large, stiff system of ordinary 
differential equations (ODEs).
The application of common implicit temporal discretization methods, such as the Crank--Nicolson scheme, requires the solution 
of very large, dense linear systems in every time step.
The main aim of the present paper is to study a variety of recent and novel operator splitting schemes that 
effectively overcome this.

For one-dimensional option valuation PIDEs, Cont \& Voltchkova \cite{CV05} proposed an {\it implicit-explicit} (IMEX) Euler 
time-stepping scheme, where the PDE part is treated implicitly and the integral part explicitly. 
This IMEX Euler scheme is only first-order consistent. 
Various authors subsequently considered higher-order IMEX schemes for PIDEs in finance, see for example \cite{B07,FL08,KL11}.
Salmi \& Toivanen \cite{ST14} recommended the IMEX CNAB scheme, where the PDE part is handled in a Crank--Nicolson manner 
and the integral part is treated in a two-step Adams--Bashforth fashion.
This second-order two-step IMEX scheme has next been applied to two-dimensional option valuation PIDEs in for example 
\cite{HT16,STS14}. 

An alternative approach was proposed by Tavella \& Randall \cite{TR00book} and d'Halluin, Forsyth \& Vetzal \cite{HFV05}.
These authors considered the Crank--Nicolson scheme together with a fixed-point iteration on the integral part for the 
solution of the linear system in each time step. 
This approach has been applied to two-dimensional option valuation PIDEs in Clift \& Forsyth \cite{CF08}, including the 
two-dimensional Merton PIDE.
It is found that in general two to three iterations are sufficient to reach a satisfactory error level.
When the number of fixed-point iterations is frozen, a particular one-step IMEX scheme is obtained, see Section~\ref{SecTime}
below.

For multidimensional option valuation PDEs (without integral part), operator splitting schemes are widely considered 
in the literature.
Notably, {\it alternating direction implicit} (ADI) methods, which combine implicit unidirectional correction stages 
with explicit stages involving the mixed derivative terms, are highly effective, see for example 
\cite{GHHP18,H17book,HT16,HW09,S19}.
Hundsdorfer \& in 't Hout \cite{HH18}~recently studied a novel class of multistep stabilizing correction methods for 
PDEs with mixed derivative terms and showed that such methods can be competitive with the well-established one-step 
modified Craig--Sneyd (MCS) scheme \cite{HW09} in the application to the familiar two-dimensional Heston PDE~\cite{H93}.

In the past years, ADI methods have been extended and applied to two-dimensional PIDEs in finance. 
In \cite{I11} Strang splitting was employed between the PDE and integral parts and an ADI scheme was used for
the PDE subproblems, see also \cite{I17}.
The direct extension of ADI methods to PIDEs has recently been considered in \cite{HT18,KLR18}.
Here the implicit unidirectional stages are combined with explicit stages involving both the mixed derivative term 
and the integral part.
In 't Hout \& Toivanen \cite{HT18} defined three adaptations of the MCS scheme to two-dimensional PIDEs that are all 
second-order consistent.  
They found that the adaptation where the integral part is treated in a two-step Adams--Bashforth fashion is 
preferred in the application to the well-known two-dimensional Bates PIDE \cite{B96}, which merges the Heston PDE and 
the one-dimensional Merton PIDE.

An outline of the rest of this paper is as follows.
Section~\ref{SecModel} describes the asset price dynamics under the two-asset Merton jump-diffusion model and introduces 
the two-dimensional PIDE that holds for the value of a rainbow option under this model.
Section~\ref{SecSpatial} deals with the first step in the MOL approach, the spatial discretization of the two-dimensional 
Merton PIDE.
Here particular attention is given to a suitable discretization of the double integral part, which enables application 
of the FFT algorithm from \cite{BTK01}.
Section~\ref{SecTime} concerns the second step in the MOL approach, the temporal discretization of the semidiscretized PIDE.
We formulate seven recent and novel operator splitting schemes of both the IMEX and the ADI kind, which all treat the 
integral part in an explicit fashion.
In Section~\ref{SecResults} extensive numerical experiments are presented.
Here two types of rainbow options are considered: a European put-on-the-min option and a European put-on-the-average option.
For all seven splitting schemes, the behaviour of the temporal discretization error is investigated in detail. 
Moreover, employing the recent semi-closed analytical formula derived in Boen \cite{B18}, also the total discretization 
error is studied in the case of the put-on-the-min option.
The final Section~\ref{SecConc} gives conclusions.


\section{The two-asset Merton jump-diffusion model}\label{SecModel}
Let $(\Omega, \mathcal{F}, \{\mathcal{F}_\tau\}, \bQ)$ be a probability measure space, with $\bQ$ a risk-neutral measure. 
Under a two-asset jump-diffusion model, the asset prices under $\bQ$ are given by
\begin{align*}
S^{(1)}_\tau &= S^{(1)}_0\exp\left((r-\tfrac{1}{2}\sigma_1^2-\lambda\kappa_1)\tau+\sigma_1W^{(1)}_\tau+\sum_{k=1}^{N_\tau}Y^{(1)}_k\right),\\
S^{(2)}_\tau &= S^{(2)}_0\exp\left((r-\tfrac{1}{2}\sigma_2^2-\lambda\kappa_2)\tau+\sigma_2W^{(2)}_\tau+\sum_{k=1}^{N_\tau}Y^{(2)}_k\right),
\end{align*}
where 
\begin{itemize}
\item $r$ is the risk-free interest rate,
\item $\sigma_i$ ($i=1,2$) is the instantaneous volatility of asset $i$, conditional on the event that no jumps occur,
\item $W_\tau = (W^{(1)}_\tau,W^{(2)}_\tau)^\top$ consists of two correlated standard Brownian motions, with
correlation coefficient~$\rho$, 
\item $\lambda$ is the jump intensity of the Poisson arrival process $N = \{N_\tau, \tau\geq0\}$, 
\item $Y^{(i)}_k (k=1,2,3,\ldots)$ are independent jump sizes that have an identical distribution to the random jump size $Y^{(i)}$
and, for any given $k$, the jumps $Y^{(i)}_k \ (i=1,2)$ occur simultaneously, driven by the same Poisson arrival process $N$, 
and are correlated with correlation coefficient~$\widehat{\rho}$,
\item $\kappa_i$ is the expected relative jump size
\[
\kappa_i = \bE_{\bQ}[e^{Y^{(i)}} -1] \quad (i=1,2).
\]
\end{itemize}
The $W_\tau$, $N_\tau$ and $Y = (Y^{(1)},Y^{(2)})^\top$ are assumed to be independent of each other whenever $\tau\ge 0$.

In the two-dimensional version of the Merton model \cite{M76}, the random vector $Y$ is bivariate normally distributed with 
mean $(\gamma_1, \gamma_2)^\top$ and covariance matrix
\[
\Sigma_{Y}=\begin{pmatrix}
\delta_1^2 & \widehat{\rho}\delta_1\delta_2\\
\widehat{\rho}\delta_1\delta_2 & \delta_2^2
\end{pmatrix},
\]
where $\delta_i^2$ is the variance of $Y^{(i)}$ $(i=1,2)$. 
The moment generating function of $Y$ is then given by
\begin{equation*}
M_{\bQ,Y}(u)  = \bE_{\bQ}\left[e^{u^\top Y}\right] = 
\exp\left(\gamma_1u_1 + \gamma_2u_2 + \tfrac{1}{2} \left(\delta_1^2u_1^2 + 
2\widehat{\rho}\delta_1\delta_2u_1u_2 + \delta_2^2u_2^2\right)\right),
\quad u = \begin{pmatrix} u_1 \\ u_2 \end{pmatrix},
\end{equation*}
and the expected relative jump size is
\[
\kappa_i= e^{\gamma_i+\half \delta_i^2}-1 \quad (i=1,2). 
\]
Under the two-asset Merton jump-diffusion model, the value $v = v(s_1,s_2,t)$ of a European-style\footnote{A European-style 
option can only be exercised at the maturity date $T$.} option with maturity date $T>0$ and $s_i$ ($i=1,2$) representing the 
price of asset $i$ at time $\tau = T-t$, satisfies the following PIDE~\cite{CF08}:
\begin{align}\label{PIDE2D}
\dud{t} =&~ \tfrac{1}{2} \sigma_1^2s_1^2\dudd{s_1} + \rho \sigma_1\sigma_2s_1s_2\duddm{s_1}{s_2} + \tfrac{1}{2} \sigma_2^2 s_2^2\dudd{s_2} + 
(r-\lambda \kappa_1) s_1\dud{s_1} + (r-\lambda\kappa_2) s_2 \dud{s_2} \nonumber \\
 & -(r+\lambda)v+\lambda\int_0^{\infty}\int_0^{\infty} v(s_1y_1, s_2y_2,t) f(y_1,y_2) \dd y_1 \dd y_2
\end{align}
whenever $s_1>0$, $s_2>0$, $0<t\leq T$.
Here
\begin{align*}
f(y_1,y_2) = \frac{1}{2\pi \delta_1\delta_2\sqrt{1-\rhoh^{\hspace{0.01cm} 2}}\,y_1y_2}\exp\left(-\frac{\left(\frac{\ln(y_1)-\gamma_1}{\delta_1}\right)^2+\left(\frac{\ln(y_2)-\gamma_2}{\delta_2}\right)^2-2\rhoh\left(\frac{\ln(y_1)-\gamma_1}{\delta_1}\right)\left(\frac{\ln(y_2)-\gamma_2}{\delta_2}\right)}{2(1-\rhoh^{\hspace{0.01cm} 2})}\right)
\end{align*}
is the probability density function of a bivariate lognormal distribution. 
The initial condition for \eqref{PIDE2D} is given by
\[
v(s_1,s_2,0) = \phi(s_1,s_2),
\]
where $\phi$ denotes the payoff function of the option, that is, the option value at expiry.
Further, as the boundary condition, it holds that the PIDE \eqref{PIDE2D} is itself fulfilled on the two sides $s_{1} = 0$ and 
$s_{2} = 0$, respectively.


\section{Spatial discretization}\label{SecSpatial}
For the numerical solution of the initial-boundary value problem for \eqref{PIDE2D} we consider the well-known and versatile method 
of lines (MOL) approach, which consists of two consecutive steps~\cite{HV03book}: first the PIDE \eqref{PIDE2D} is discretized in 
space and subsequently in time.
This section discusses the spatial discretization. In the next section we turn to the temporal discretization.


\subsection{Convection-diffusion-reaction part}
To render the numerical solution of \eqref{PIDE2D} feasible, the spatial domain is truncated to a bounded set 
$[0,S_{\rm max}]\times[0,S_{\rm max}]$ with fixed value $S_{\rm max}$ taken sufficiently large.
On the two far sides $s_{1} = S_{\rm max}$ and $s_{2} = S_{\rm max}$ a linear boundary condition is imposed,
which is common in finance,
\begin{equation}\label{LBC}
\dudd{s_1} = 0~~(\textrm{if}~s_{1} = S_{\rm max}) \quad \mbox{ and } \quad  \dudd{s_2} = 0 ~~(\textrm{if}~s_{2} = S_{\rm max}).
\end{equation}

In this subsection we describe the finite difference discretization of the convection-diffusion-reaction part
of \eqref{PIDE2D}, given by
\begin{equation*}\label{lindiffop}
\mathcal{D} v =  \tfrac{1}{2} \sigma_1^2s_1^2\dudd{s_1} + \rho \sigma_1\sigma_2s_1s_2\duddm{s_1}{s_2} + \tfrac{1}{2} \sigma_2^2 s_2^2\dudd{s_2} 
+ (r-\lambda \kappa_1) s_1\dud{s_1} + (r-\lambda\kappa_2) s_2 \dud{s_2} -(r+\lambda)v.
\end{equation*}
The option value function $v$ will be approximated at a nonuniform, Cartesian set of spatial grid points,
\begin{equation}\label{sgrid}
(s_{1,i},s_{2,j})\in [0,S_{\rm max}]\times[0,S_{\rm max}] \quad (0\leq i \leq m_1, \ 0\leq j \leq m_2),
\end{equation}
with $s_{1,0} = s_{2,0} = 0$ and $s_{1,m_1} = s_{2,m_2} = S_{\rm max}$.
The nonuniform grid in each spatial direction is chosen such that relatively many grid points lie in a region 
of financial and/or numerical interest.
It is generated by applying a smooth transformation to an artificial uniform grid. 
In Section~\ref{SecResults} more details shall be given of the specific type of spatial grids that is used.

Let $0=s_0<s_1<\cdots<s_m=S_{\rm max}$ be any given smooth, nonuniform, unidirectional spatial grid with mesh widths
$h_i = s_i-s_{i-1}$ ($1\leq i \leq m$) and let $u: [0,S_{\rm max}] \to \mathbb{R}$ be any given smooth function. 
The following second-order central finite difference formulas are employed for the approximation of the first
and second derivatives of $u$:
\begin{align*}\label{semconv}
u^\prime (s_i) \approx \omega_{i,-1} u(s_{i-1}) + \omega_{i,0} u(s_{i}) + \omega_{i,1} u(s_{i+1}),
\end{align*}
with 
\[
\omega_{i,-1} = \frac{-h_{i+1}}{h_{i}(h_{i}+h_{i+1})},\quad
\omega_{i,0} = \frac{h_{i+1}-h_{i}}{h_{i}h_{i+1}},\quad 
\omega_{i,1} = \frac{h_{i}}{h_{i+1}(h_{i}+h_{i+1})}
\]
and
\begin{align*}\label{semdiff}
u^{\prime\prime} (s_i) \approx \omega_{i,-1} u(s_{i-1}) + \omega_{i,0} u(s_{i}) + \omega_{i,1} u(s_{i+1}),
\end{align*}
with
\[
\omega_{i,-1} = \frac{2}{h_{i}(h_{i}+h_{i+1})},\quad 
\omega_{i,0} = \frac{-2}{h_{i}h_{i+1}},\quad 
\omega_{i,1} = \frac{2}{h_{i+1}(h_{i}+h_{i+1})}.
\]
for $1\le i \le m-1$.
If $i=0$, no finite difference formulas are required in view of the degeneracy of \eqref{PIDE2D} at the zero boundaries.
If $i=m$, the first derivative is approximated by the first-order backward finite difference formula and the second derivative
vanishes by the linear boundary condition \eqref{LBC}.
Concerning the mixed derivative term in \eqref{PIDE2D}, this is discretized by successively applying the relevant finite 
difference formulas for the first derivative in the two spatial directions.

Let the vector $V(t) = (V_{0,0}(t), V_{1,0}(t),\ldots,V_{m_1-1,m_2}(t),V_{m_1,m_2}(t))^{\top}$ where entry $V_{i,j}(t)$ 
denotes the semidiscrete approximation to $v(s_{1,i},s_{2,j},t)$ whenever $0\leq i \leq m_1$, $0\leq j \leq m_2$.
The semidiscrete version of the convection-diffusion-reaction part $\mathcal{D} v$ can then be written as 
\[
A^{(D)} V(t)
\] 
with matrix
\[
A^{(D)} = A^{(M)} + A_1 + A_2,
\]
where
\begin{align*}
    A^{(M)} &= \rho\sigma_1\sigma_2\left(X_2D_2^{(1)}\right)\otimes\left(X_1D_1^{(1)}\right),\\
    A_1 &= I_2 \otimes \left(\tfrac{1}{2}\sigma_1^2X_1^2D_1^{(2)} + (r-\lambda\kappa_1)X_1D_1^{(1)}-\tfrac{1}{2}(r+\lambda)I_1\right),\\
    A_2 &= \left(\tfrac{1}{2}\sigma_2^2X_2^2D_2^{(2)} + (r-\lambda\kappa_2)X_2D_2^{(1)}-\tfrac{1}{2}(r+\lambda)I_2\right)\otimes I_1.
\end{align*}
Here $I_k$, $X_k$, $D_k^{(l)}$ are given $(m_k+1)\times (m_k+1)$ matrices for $k,l\in \{1,2\}$, where $I_k$ is the identity matrix, 
$X_k$ is the diagonal matrix
\begin{equation*}
X_k = {\rm diag}(s_{k,0}, s_{k,1},\ldots, s_{k,m_k})
\end{equation*}
and $D_k^{(l)}$ is the matrix representing numerical differentiation of order $l$ in the $k$-th spatial direction by the relevant 
finite difference formula above.
The matrix $A^{(M)}$ corresponds to the mixed derivative term in \eqref{PIDE2D} and $A_k$ corresponds to all derivative terms in 
the $k$-th spatial direction ($k=1,2$). 
Further, the reaction term has been distributed equally across $A_1$ and $A_2$.


\subsection{Integral part}
In order to discretize the integral part in the PIDE \eqref{PIDE2D}, we apply a well-known transformation to the 
log-price variable. 
Let $x_i = \ln(s_i)$ be the log-price of the $i$-th asset and $\eta_i = \ln(y_i)$ ($i=1,2$), then 
\begin{align*}
    \mathcal{J}(s_1,s_2):=&\ \lambda\int_0^{\infty}\int_0^{\infty}v(s_1y_1,s_2y_2,t)f(y_1,y_2) \dd y_1 \dd y_2\\
    =&\ \lambda\int_{-\infty}^{\infty}\int_{-\infty}^{\infty} \overline{v}(x_1+\eta_1,x_2+\eta_2,t) \overline{f}(\eta_1,\eta_2)\dd \eta_1 \dd \eta_2\\
    =&\ \lambda\int_{-\infty}^{\infty}\int_{-\infty}^{\infty} \overline{v}(\xi_1,\xi_2,t) \overline{f}(\xi_1-x_1,\xi_2-x_2)\dd \xi_1 \dd \xi_2,
\end{align*}
where $\overline{v}(\xi_1,\xi_2,t)=v(e^{\xi_1},e^{\xi_2},t)$ and $\overline{f}(\eta_1,\eta_2) =$  $f(e^{\eta_1},e^{\eta_2})e^{\eta_1}e^{\eta_2}$ is the probability density function of a bivariate normal distribution. Since $\overline{f}$ is a real-valued function, it is clear that the double integral in the log-price variable can be viewed as a two-dimensional cross-correlation (see e.g.~\cite{B2000book}):
\[
\overline{\mathcal{J}}(x_1,x_2) := 
\lambda\int_{-\infty}^{\infty}\int_{-\infty}^{\infty} \overline{v}(\xi_1,\xi_2,t) \overline{f}(\xi_1-x_1,\xi_2-x_2)\dd \xi_1 \dd \xi_2 =
\lambda (\overline{v} \star \star \overline{f})(x_1,x_2),
\]
where $\star\star$ denotes the two-dimensional cross-correlation operator and $\overline{\mathcal{J}}(x_1,x_2) = \mathcal{J}(e^{x_1},e^{x_2})$. For the semidiscretization, we restrict the log-price domain to $[-X_{\rm max},X_{\rm max}]\times [-X_{\rm max},X_{\rm max}]$ with~$X_{\rm max} = \ln(S_{\rm max})$ and consider a $M_1 \times M_2$ uniform grid $(x_{1,k},x_{2,l}) = (k\Delta x_1,l\Delta x_2)$ for ${k=-\widecheck{M}_1+1,\ldots,\widecheck{M}_1, \ l=-\widecheck{M}_2+1,\ldots,\widecheck{M}_2}$, where $\widecheck{M}_i = M_i/2$, 
$\widecheck{M}_i\Delta x_i = X_{\rm max}$ and $M_i$ is chosen to be a power of 2 such that the mesh width $\Delta x_i$ is smaller than the smallest mesh width in the nonuniform $\ln(s_i)$-grid ($i=1,2$). The double integral $\overline{\mathcal{J}}(x_{1,k},x_{2,l})$ is approximated by
\[
\overline{J}_{k,l} =\lambda \sum_{j=-\widecheck{M}_2+1}^{\widecheck{M}_2}\sum_{i=-\widecheck{M}_1+1}^{\widecheck{M}_1} \overline{V}_{i,j}\overline{F}_{i-k,j-l},
\]
with approximation $\overline{V}_{i,j} \approx \oV(x_{1,i}, x_{2,j},t)$ defined below and $\overline{F}_{i-k,j-l}=\overline{f}_{i-k,j-l} \Delta x_1 \Delta x_2$ where
\begin{align*}
\overline{f}_{i-k,j-l} &= \overline{f}(x_{1,i}-x_{1,k},x_{2,j}-x_{2,l}).
\end{align*}
Let $E(M,d)$ denote the $M\times M$ matrix with ones on the $d$-th off-diagonal and zeros elsewhere. The matrix 
\[
\overline{A}^{(J)} = \lambda \overline{F}
\]
with 
\[
~~~\overline{F} = \sum_{d=-M_2+1}^{M_2-1} E(M_2, d) \otimes \overline{F}^1_{d}
\]
and
\[
\overline{F}^1_{d} =  \sum_{c=-M_1+1}^{M_1-1} E(M_1, c) \overline{F}_{c,d}
\]
is then an asymmetrical two-level block-Toeplitz matrix of size $ (M_1M_2) \times (M_1M_2)$. Its structure is illustrated in Appendix \ref{appToep}.
Let the vectors
\begin{align*}
\overline{J}=&\ (\overline{J}_{-\widecheck{M}_1+1, -\widecheck{M}_2+1}, \overline{J}_{-\widecheck{M}_1+2, -\widecheck{M}_2+1}, \ldots, \overline{J}_{\widecheck{M}_1-1, \widecheck{M}_2}, \overline{J}_{\widecheck{M}_1, \widecheck{M}_2})^\top,\\\\
\overline{V} =&\ (\overline{V}_{-\widecheck{M}_1+1, -\widecheck{M}_2+1}, \overline{V}_{-\widecheck{M}_1+2, -\widecheck{M}_2+1}, \ldots, \overline{V}_{\widecheck{M}_1-1, \widecheck{M}_2}, \overline{V}_{\widecheck{M}_1, \widecheck{M}_2})^\top,
\end{align*}
then
\[
\overline{J} =\overline{A}^{(J)} \ \overline{V}.
\]
This matrix-vector product can be computed very efficiently using the FFT algorithm for matrix-vector multiplications with asymmetric multilevel block-Toeplitz matrices, introduced in Barrowes, Teixeira \& Kong \cite{BTK01}. This fast algorithm computes the matrix-vector product using two FFTs and one inverse FFT (IFFT), reducing the computational cost to $\mathcal{O}(M_1M_2\log(M_1M_2))$. The algorithm from \cite{BTK01} essentially embeds the block-Toeplitz matrix in a circulant matrix before applying the FFT, conveniently avoiding any wrap-around effect, which would have been present when applying the FFT directly to the block-Toeplitz matrix. Further, we mention that the algorithm from \cite{BTK01} is more efficient than that employed in \cite{STS14}, as the former uses only one-dimensional FFTs, whereas the latter requires two-dimensional FFTs.

Note that the spatial grid on which we discretize the PIDE \eqref{PIDE2D} is nonuniform, and in general does not lead to a uniform log-price grid. We therefore (bi)linearly interpolate the option value approximations between the two grids directly before and after application of the algorithm of \cite{BTK01}. Linearly interpolating $V$ to $\overline{V}$ (using Lagrange basis functions on the $s$-grid, as in e.g. \cite{HFV05}), yields
\begin{align}\label{interpV}
\overline{V}_{i,j} = \ &\phi_{\Pi(i)}\psi_{\Omega(j)}V_{\Pi(i),\Omega(j)} + (1- \phi_{\Pi(i)}) \psi_{\Omega(j)}V_{\Pi(i)+1,\Omega(j)}\nonumber \\ &+  \phi_{\Pi(i)}(1-\psi_{\Omega(j)})V_{\Pi(i),\Omega(j)+1} +  (1-\phi_{\Pi(i)})(1-\psi_{\Omega(j)})V_{\Pi(i)+1,\Omega(j)+1},
\end{align}
whenever $s_{1,\Pi(i)} \leq e^{x_{1,i}} \leq s_{1,\Pi(i)+1}$ and $s_{2,\Omega(j)} \leq e^{x_{2,j}} \leq s_{1,\Omega(j)+1}$, where $\phi_{\Pi(i)},\psi_{\Omega(j)} \in [0,1]$ are interpolation weights . This can be written as
\[
\overline{V} = \mathcal{X} V,
\]
where the matrix $\mathcal{X}$ contains the interpolation weights such that \eqref{interpV} holds. Similarly, linearly interpolating $\overline{J}$ back to the original $s$-grid, leads to approximations $J_{k,l}$ of $\mathcal{J}(s_{1,k},s_{2,l})$ given by
\begin{align}\label{interpJ}
 J_{k,l} =& \overline{\phi}_{\Lambda(k)}\overline{\psi}_{\Upsilon(l)} \overline{J}_{\Lambda(k),\Upsilon(l)} +  (1-\overline{\phi}_{\Lambda(k)})\overline{\psi}_{\Upsilon(l)} \overline{J}_{\Lambda(k)+1,\Upsilon(l)}\nonumber \\
 &+  \overline{\phi}_{\Lambda(k)}(1-\overline{\psi}_{\Upsilon(l)}) \overline{J}_{\Lambda(k),\Upsilon(l)+1} +  (1-\overline{\phi}_{\Lambda(k)})(1-\overline{\psi}_{\Upsilon(l)}) \overline{J}_{\Lambda(k)+1,\Upsilon(l)+1},  
\end{align}
whenever $e^{x_{1,\Lambda(k)}} \leq s_{1,k} \leq  e^{x_{1,\Lambda(k)+1}}$ and $e^{x_{2,\Upsilon(l)}} \leq s_{2,l} \leq  e^{x_{2,\Upsilon(l)+1}}$, where $\overline{\phi}_{\Lambda(k)},\overline{\psi}_{\Upsilon(l)}  \in [0,1]$ are interpolation weights. This can be written as
\[
J = \overline{\mathcal{X}} \ \overline{J},
\]
where the matrix $\overline{\mathcal{X}}$ contains the interpolation weights such that \eqref{interpJ} holds. Consequently, we have
\[
J =  \overline{\mathcal{X}} \ \overline{J} = \overline{\mathcal{X}}\ \overline{A}^{(J)} \overline{V} =  \overline{\mathcal{X}}\ \overline{A}^{(J)} \mathcal{X} V.
\]
The entire procedure is summarized in Algorithm \ref{algo1}.

\begin{algorithm}
\caption{Computing the approximation to the integral part in \eqref{PIDE2D}}
\label{algo1}
\begin{algorithmic}
\STATE \textbf{1.} Linearly interpolate $V(t) = (V_{0,0}(t),V_{1,0}(t),\ldots, V_{m_1-1,m_2}(t),V_{m_1,m_2}(t))^{\top}$ onto the uniform log-price grid, leading to the vector $\overline{V}$.
\STATE \textbf{2.} Extend the log-price grid $(x_{1,k},x_{2,l})$ for ${k=-\widecheck{M}_1+1,\ldots,\widecheck{M}_1,} \ {l=-\widecheck{M}_2+1,\ldots,\widecheck{M}_2}$ on both sides, keeping it uniform:
\[
x_{1,-M_1+1},\ldots, x_{1,-\widecheck{M}_1}, \ldots, x_{1,0},\ldots,x_{1,\widecheck{M}_1},\ldots, x_{1,M_1-1},
\]
\[
x_{2,-M_2+1},\ldots, x_{2,-\widecheck{M}_2}, \ldots, x_{2,0},\ldots,x_{2,\widecheck{M}_2},\ldots, x_{2,M_2-1},
\]
and compute $\of_{i-k,j-l}$ for $i,k=-\widecheck{M}_1+1,\ldots,\widecheck{M}_1, \ j,l=-\widecheck{M}_2+1,\ldots,\widecheck{M}_2$.
\STATE \textbf{3.} Compute the matrix-vector product $\overline{A}^{(J)} \overline{V}$ using the fast algorithm of \cite{BTK01}, resulting in the vector~$\overline{J}$.
\STATE \textbf{4.} Linearly interpolate $\overline{J}$ back to the original nonuniform grid $(s_{1,i},s_{2,j})$ for $i=1,\ldots m_1$, $j=1,\ldots m_2$, leading to the vector $J$.
\end{algorithmic}
\end{algorithm}

Note that in the interpolation step from the uniform log-price grid to the nonuniform price grid, we obtain values $J_{i,j}$ for $i=1,\ldots,m_1, \ j=1,\ldots,m_2$. On the boundaries $s_1 = 0$ and $s_2 = 0$, the two-dimensional Merton PIDE \eqref{PIDE2D} reduces to the one-dimensional Merton PIDE and the jump part can be computed using a similar FFT technique (see f.i.~\cite{HFV05}).

Denoting the approximation of $\mathcal{J}$ on the full grid \eqref{sgrid} by $A^{(J)}V(t)$, we arrive at the 
following semidiscrete system of ODEs:
\begin{align}\label{ODEs}
V'(t) = AV(t),
\end{align}
for $0<t\leq T$, where $A = A^{(D)} + A^{(J)} = A^{(M)} + A_1 + A_2 + A^{(J)}$.

The maximum norm of the matrix $A^{(J)}$ can be bounded independently of the spatial grid by a moderate constant whenever $\lambda$ 
is moderate.
Accordingly, the dense matrix $A^{(J)}$ constitutes a nonstiff part of the semidiscrete system, whereas the sparse matrix $A^{(D)}$ 
constitutes a stiff part.


\section{Temporal discretization}\label{SecTime}
For the temporal discretization of the semidiscrete system \eqref{ODEs} we shall investigate seven operator splitting schemes.
Each of these schemes conveniently treats the nonstiff integral part, given by the matrix $A^{(J)}$, in an explicit fashion.

Let integer $N\ge 1$ be given and step size $\Delta t = T/N$.
Each of the following schemes defines an approximation $V^{n}$ to $V(t_n)$ at the temporal grid point $t_n = n \Delta t$ 
successively for $n=1,2,\ldots,N$.
The initial vector $V(0)=V^0$ is determined by the payoff function and will be specified in Section~\ref{SecResults}.

\begin{enumerate}

    \item \textit{Crank--Nicolson Forward Euler (CNFE) scheme}:
    \vskip0.1cm
    In this basic method the convection-diffusion-reaction part, including the mixed derivative, is handled implicitly 
    using the Crank--Nicolson scheme and the integral part is treated explicitly in a forward Euler manner:
    \begin{align}\label{CNFE}
    \left(I-\tfrac{1}{2}\Delta t A^{(D)}\right)V^{n} = \left(I + \tfrac{1}{2}\Delta t A^{(D)}\right)V^{n-1} + \Delta t A^{(J)} V^{n-1}.
    \end{align}
    Method \eqref{CNFE} is a one-step IMEX scheme.
    Due to the application of forward Euler, its order\footnote{This refers to the classical order of consistency, that 
    is, for fixed nonstiff ODEs.} is just equal to one. 
    \vskip0.2cm
    \item \textit{Crank--Nicolson scheme with fixed-point iteration (CNFI)}:
    \vskip0.1cm
    In Tavella \& Randall \cite{TR00book} and d'Halluin, Forsyth \& Vetzal \cite{HFV05} the Crank--Nicolson scheme has been 
    proposed with a fixed-point iteration on the integral part:
    \begin{align}\label{CNFI}
    \left(I-\tfrac{1}{2}\Delta t A^{(D)}\right)Y_k = \left(I + \tfrac{1}{2}\Delta t A^{(D)}\right)V^{n-1} +  \tfrac{1}{2} \Delta t A^{(J)} (Y_{k-1}+V^{n-1})
    \end{align}
    for $k=1,2,\ldots,l$ and $V^n = Y_l$ using the starting value $Y_0 = V^{n-1}$. 
    Clearly, with just one iteration, method \eqref{CNFE} is obtained.
    Method \eqref{CNFI} can be applied with a dynamic convergence criterion. Numerical experiments in \cite{CF08,HFV05}
    reveal that this generally leads to $l=2$ or $l=3$ iterations.
    For comparable computational work to the other schemes, formulated below, we choose here a fixed number of $l=2$ 
    iterations.
    Then method \eqref{CNFI} can be regarded as a one-step IMEX scheme, where the integral part is treated by the 
    explicit trapezoidal rule, also called the modified Euler method. 
    It can be verified, by Taylor expansion, that its order is equal to two.
    \vskip0.2cm
    \item \textit{Implicit-Explicit Trapezoidal Rule (IETR)}:
    \vskip0.1cm
    A different combination of the implicit trapezoidal rule (Crank--Nicolson) for the convection-diffusion-reaction part
    and its explicit variant for the integral part is considered by in 't Hout \cite{H17book}:
    \begin{align}\label{IETR}
    \begin{cases}
    Y_0 = V^{n-1} + \Delta t\, (A^{(D)}+A^{(J)})V^{n-1},\\
    \widehat{Y}_0 = Y_0+\tfrac{1}{2} \Delta t A^{(J)} \big(Y_0-V^{n-1}\big),\\
    Y_1 = \widehat{Y}_0+\tfrac{1}{2} \Delta t A^{(D)} \big(Y_1-V^{n-1}\big),\\
    V^n = Y_1.
    \end{cases}
    \end{align}
    Here the first two internal stages $Y_0$, $\widehat{Y}_0$ are explicit, whereas the third internal stage $Y_1$ is 
    implicit.
    Method \eqref{IETR} is also a one-step IMEX scheme with order equal to two.
    This method is already well-known in the literature on the numerical solution of PDEs, without integral part, 
    see e.g. Hundsdorfer \& Verwer \cite{HV03book}.
    \vskip0.2cm    
    \item \textit{Crank--Nicolson Adams--Bashforth (CNAB) scheme}:
    \vskip0.1cm
    In the CNAB method, the convection-diffusion-reaction part is again treated by the Crank--Nicolson scheme, but the 
    integral part is now handled in a two-step Adams--Bashforth fashion:
    \begin{align}\label{CNAB}
        \left(I-\tfrac{1}{2} \Delta t A^{(D)}\right)V^{n} = \left(I+\tfrac{1}{2} \Delta t A^{(D)}\right)V^{n-1} + \tfrac{1}{2}\Delta t A^{(J)}(3V^{n-1}-V^{n-2}).
    \end{align}
    Method \eqref{CNAB} forms a two-step IMEX scheme and is of order equal to two.
    This method has first been studied for the application to PIDEs in Salmi \& Toivanen \cite{ST14} and 
    Salmi, Toivanen \& von Sydow \cite{STS14}. 
    Previously, it has been investigated for the numerical solution of PDEs in e.g.~\cite{FHV97,HV03book,H02}.
    \vskip0.2cm 
\end{enumerate}

It is well-known that the nonsmoothness of the initial (payoff) function can have an adverse effect on the convergence
behaviour of the Crank--Nicolson scheme, which can be resolved by first applying, at initial time $t=0$, two 
time steps with step size $\Delta t/2$ using the implicit Euler method to define the approximation $V^1$ to
$V(t_1)$, see e.g.~\cite{H17book,PVF03,R84}.
In the present application to PIDEs, the implicit Euler method can be computationally demanding.
We therefore apply, for all four schemes above, the IMEX Euler scheme as the starting method:
    \begin{align*}    
    \left(I-\tfrac{1}{2}\Delta t A^{(D)}\right)V^{\tfrac{1}{2}} &= V^{0} + \tfrac{1}{2}\Delta t A^{(J)} V^{0},\\
    \left(I-\tfrac{1}{2}\Delta t A^{(D)}\right)V^{1} &= V^{\tfrac{1}{2} } + \tfrac{1}{2}\Delta t A^{(J)} V^{\tfrac{1}{2}}.
    \end{align*} 

In the following, three recent operator splitting schemes for PIDEs are formulated that employ a subsequent, useful 
splitting of the two-dimensional convection-diffusion-reaction part, given by the matrix $A^{(D)}$.

\begin{enumerate} 
    \item[5.] \textit{One-step adaptation of the modified Craig--Sneyd (MCS) scheme}:
    \vskip0.1cm
    This method is a direct adaptation of a well-established one-step ADI scheme:
    \begin{align}\label{MCS}
    \begin{cases}
    Y_0 = V^{n-1} + \Delta t\, (A^{(D)}+A^{(J)})V^{n-1},\\
    Y_j = Y_{j-1} + \theta\Delta t A_j(Y_j-V^{n-1}) \quad (j=1,2),\\
    \widehat{Y}_0 = Y_0 + \theta \Delta t\, (A^{(M)}+A^{(J)})(Y_2 - V^{n-1}),\\
    \widetilde{Y}_0 = \widehat{Y}_0 + (\tfrac{1}{2}-\theta)\Delta t\, (A^{(D)}+A^{(J)})(Y_2-V^{n-1}),\\
    \widetilde{Y}_j = \widetilde{Y}_{j-1} + \theta \Delta t A_j(\widetilde{Y}_j - V^{n-1}) \quad (j=1,2),\\
    V^n = \widetilde{Y}_2,
    \end{cases}
    \end{align}
    where $\theta>0$ is a given parameter.
    The MCS scheme was introduced by in 't Hout \& Welfert \cite{HW09} for the numerical solution of 
    PDEs with mixed derivative terms.
    The above, direct adaptation to PIDEs has recently been studied in in 't Hout \& Toivanen
    \cite{HT18}.
    Method \eqref{MCS} is of order two for any value~$\theta$.
    Here we make the common choice $\theta = \frac{1}{3}$, which is motivated by stability and accuracy results 
    in the literature for two-dimensional problems, see e.g.~\cite{HM11,HT18,HW09,HW16}.
    It is readily seen that in \eqref{MCS} the integral part and mixed derivative term are both treated 
    by the explicit trapezoidal rule.
    Notice that the explicit stages $\widehat{Y}_0$, $\widetilde{Y}_0$ can be combined, so that the integral part 
    is evaluated twice per time step.
    The implicit stages $Y_j$, $\widetilde{Y}_j$ (for $j=1,2$) are often called stabilizing corrections and are 
    unidirectional.
    The pertinent linear systems for these stages are tridiagonal, so that they can be solved very efficiently 
    using an a priori $LU$ factorization.\\
    \vskip0.2cm    
    \item[6.] \textit{Two-step adaptation of the MCS (MCS2) scheme}:
    \vskip0.1cm
    In \cite{HT18} an alternative adaptation of the MCS scheme to PIDEs has been proposed, where the integral part 
    is dealt with in a two-step Adams--Bashforth fashion:
    \begin{align}\label{MCS2}
    \begin{cases}
    X_0 = V^{n-1}+\Delta t A^{(D)}V^{n-1},\\
    Y_0 = X_0 + \tfrac{1}{2}\Delta t A^{(J)}(3V^{n-1}-V^{n-2}),\\
    Y_j = Y_{j-1}+\theta\Delta t A_j(Y_j-V^{n-1})\quad (j=1,2),\\
    \widehat{Y}_0 = Y_0 + \theta\Delta t A^{(M)}(Y_2-V^{n-1}),\\
    \widetilde{Y}_0 = \widehat{Y}_0 + (\tfrac{1}{2}-\theta)\Delta t A^{(D)}(Y_2-V^{n-1}),\\
    \widetilde{Y}_j = \widetilde{Y}_{j-1} + \theta \Delta t A_j(\widetilde{Y}_j-V^{n-1})\quad (j=1,2),\\
    V^n = \widetilde{Y}_2.
    \end{cases}
    \end{align}
    Method \eqref{MCS2} is also of order two for any value~$\theta$.
    We take again $\theta = \frac{1}{3}$ and, for starting \eqref{MCS2}, define $V^1$ by the one-step method 
    \eqref{MCS}.
    \vskip0.2cm    
    \item[7.] \textit{Stabilizing correction two-step Adams-type (SC2A) scheme}:
    \vskip0.1cm
    In Hundsdorfer \& in 't Hout \cite{HH18} a novel class of stabilizing correction multistep methods 
    has recently been investigated for the numerical solution of PDEs.
    We select here a prominent member of this class, the two-step Adams-type scheme called SC2A, and consider its
    direct adaptation to PIDEs:
    \begin{align}\label{SC2A}
    \begin{cases}
    Y_0 = V^{n-1} + \Delta t\, (A^{(M)}+A^{(J)})\sum_{i=1}^2\widehat{b}_iV^{n-i} + \Delta t\, (A_1+A_2)\sum_{i=1}^2\widecheck{b}_iV^{n-i},\\
    Y_j = Y_{j-1} + \theta \Delta t A_j(Y_j-V^{n-1})\quad (j=1,2),\\
    V^n = Y_2,
    \end{cases}
    \end{align}
    with coefficients $(\widehat{b}_1,\widehat{b}_2) = \left(\frac{3}{2},-\half\right)$ and 
    $(\widecheck{b}_1,\widecheck{b}_2) = \left(\frac{3}{2}-\theta, - \half+\theta\right)$.
    The integral part and mixed derivative term are now both handled by the two-step Adams--Bashforth scheme.
    Method \eqref{SC2A} is also of order two for any value~$\theta$. 
    Following \cite{HH18} we take $\theta = \frac{3}{4}$, which is motivated by stability and accuracy results.
    For starting \eqref{SC2A}, the one-step method \eqref{MCS} is used with $\theta = \frac{1}{3}$ to define $V^1$.
    \vskip0.2cm    

\end{enumerate}

In each of the seven time-stepping schemes described above, any matrix-vector product with the matrix $A^{(J)}$ means 
applying Algorithm~1, hence requiring two FFTs and one IFFT. 
In the schemes \eqref{CNFE}, \eqref{CNAB}, \eqref{MCS2}, \eqref{SC2A} only one such matrix-vector product arises per time 
step, whereas two such products appear in the schemes \eqref{CNFI}, \eqref{IETR}, \eqref{MCS}.


\section{Numerical study}\label{SecResults}
In this section we examine through extensive numerical experiments the convergence behaviour of the seven operator
splitting schemes formulated in Section~\ref{SecTime} in the numerical solution of the semidiscrete two-dimensional Merton PIDE.
For the numerical study we consider the three parameter sets for the two-asset Merton jump-diffusion model and European-style option
given in Table~\ref{tabpars}.
Set~1 is taken from Clift \& Forsyth \cite{CF08}. 
Set~2 has the same diffusion parameters as in Zvan, Forsyth \& Vetzal \cite{ZFV01}, and is complemented with jump 
parameters where the intensity $\lambda$ is taken larger than in Set~1. 
Set~3 is fully new and yields more intensive jumps. 
Here $\lambda T$ is rather large, which poses a particular numerical challenge.
Notice that for all three sets the correlation coefficients $\rho$ and $\rhoh$ are both nonzero.
\begin{table}[h!]
\centering
\begin{tabular}{ c | c   c  c c c c c c c c c c }
 & $\sigma_1$ & $\sigma_2$ & $\rho$ & $\lambda$ & $\gamma_1$ & $\gamma_2$ & $\rhoh$ & $\delta_1$ & $\delta_2$ & $r$ & $K$ & $T$ \\
 \hline
Set~1 & 0.12 & 0.15 & 0.30 & 0.60 & -0.10 & 0.10 & -0.20 & 0.17 & 0.13 & 0.05 & 100 & 1\\
Set~2 & 0.30 & 0.30 & 0.50 & 2 & -0.50 & 0.30 & -0.60 & 0.40 & 0.10 & 0.05 & 40 & 0.5\\
Set~3 & 0.20 & 0.30 & 0.70 & 8 & -0.05 & -0.20 & 0.50 & 0.45 & 0.06 & 0.05 & 40 & 1
 \end{tabular}
\caption{Parameter sets for the two-asset Merton jump-diffusion model and European option.}\label{tabpars}
\end{table}
\hspace{0.5cm}\\
\\
Two types of rainbow options are considered: a European put-on-the-min option and a European put-on-the-average option. 
Their payoff functions are, respectively,
\begin{align*}
    \phi_{\text{put-on-min}}(s_1,s_2) = \max(0,K-\min(s_1,s_2))
\end{align*}
and
\begin{align*}
    \phi_{\text{put-on-average}}(s_1,s_2) = \max\left(0,K-\frac{s_1+s_2}{2}\right).
\end{align*}
 
\begin{figure}[h!]

\centering
\hspace{-0.5cm}\includegraphics[trim={1cm 0 1cm 0},clip,width=0.5\textwidth]{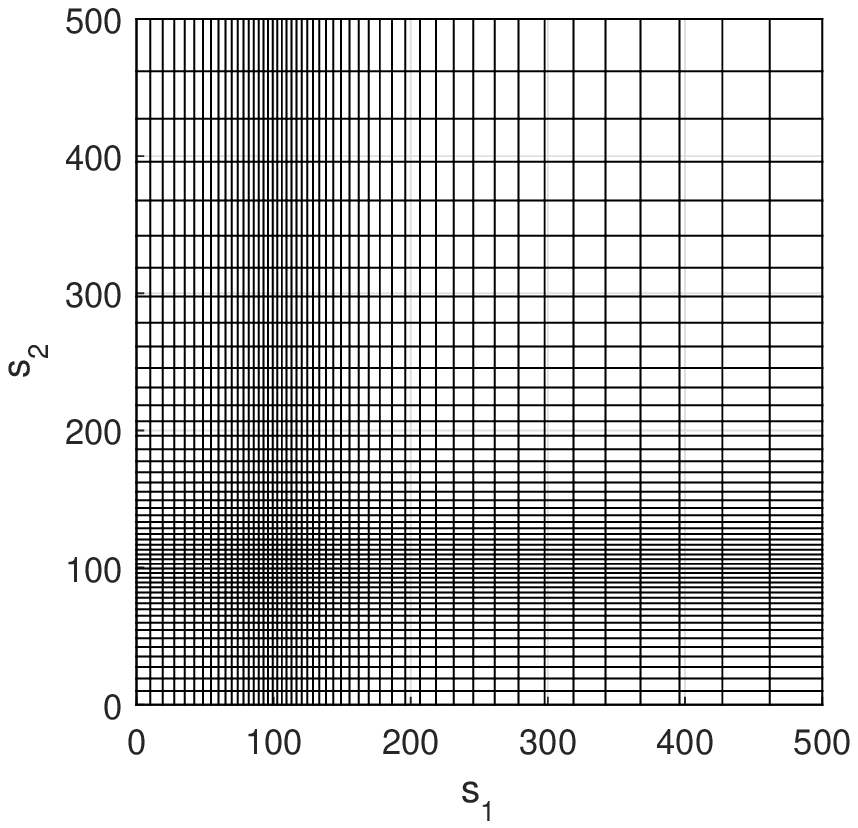}
\includegraphics[trim={1cm 0 1cm 0},clip,width=0.5\textwidth]{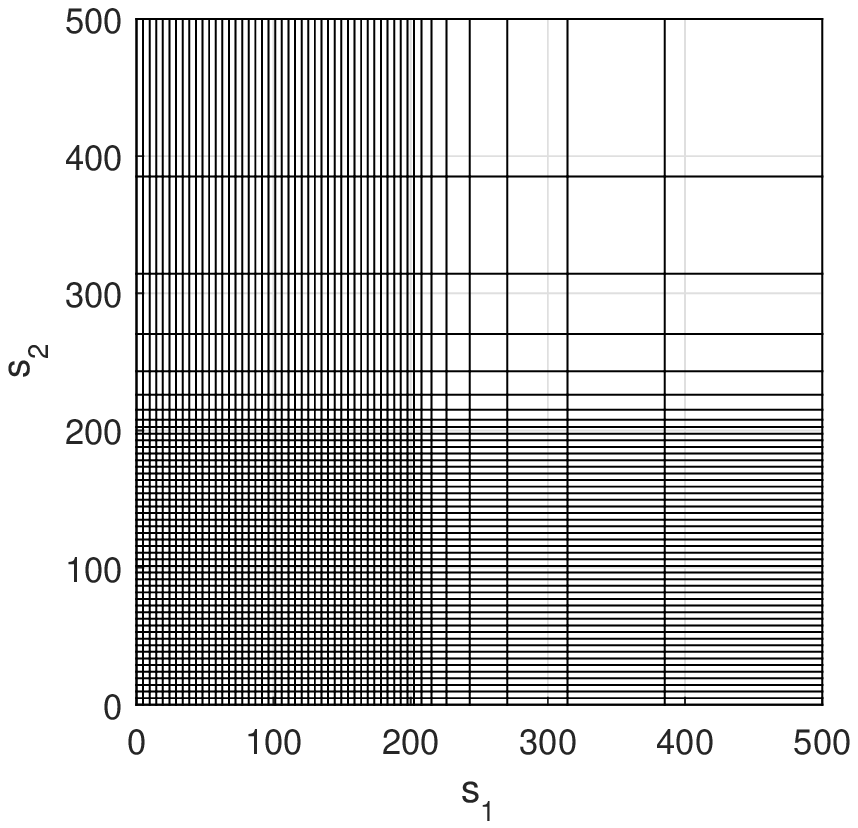}
\caption{Spatial grid for the European put-on-the-min (left) and the European put-on-the-average (right) option if $m=50$, $K = 100$, $S_{\rm max} = 5K$.}
\label{figgrid}
\end{figure}

\begin{figure}[h!]
    \centering
    \hspace{-0.5cm}\includegraphics[width =0.5\textwidth]{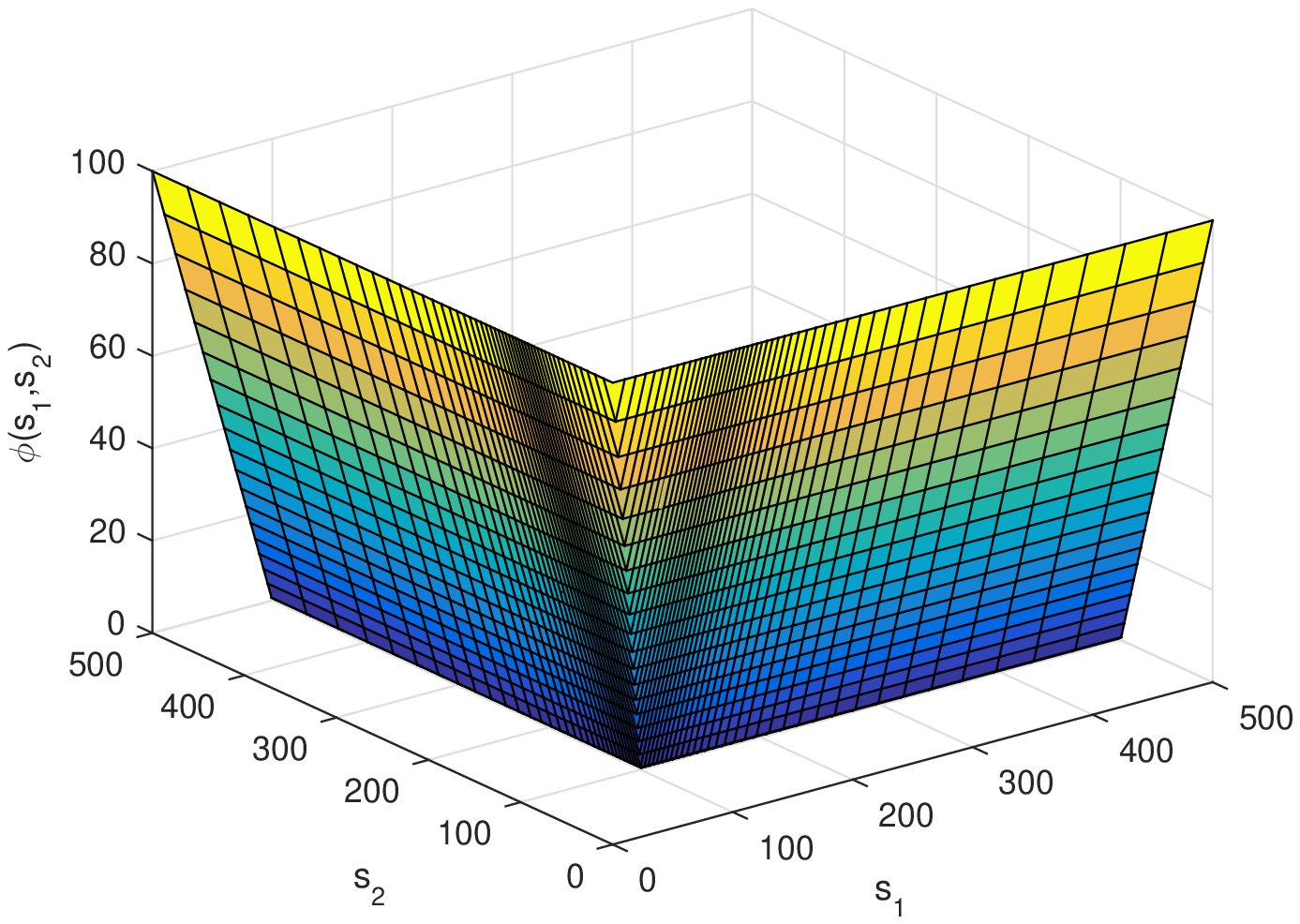}
     \includegraphics[width =0.5\textwidth]{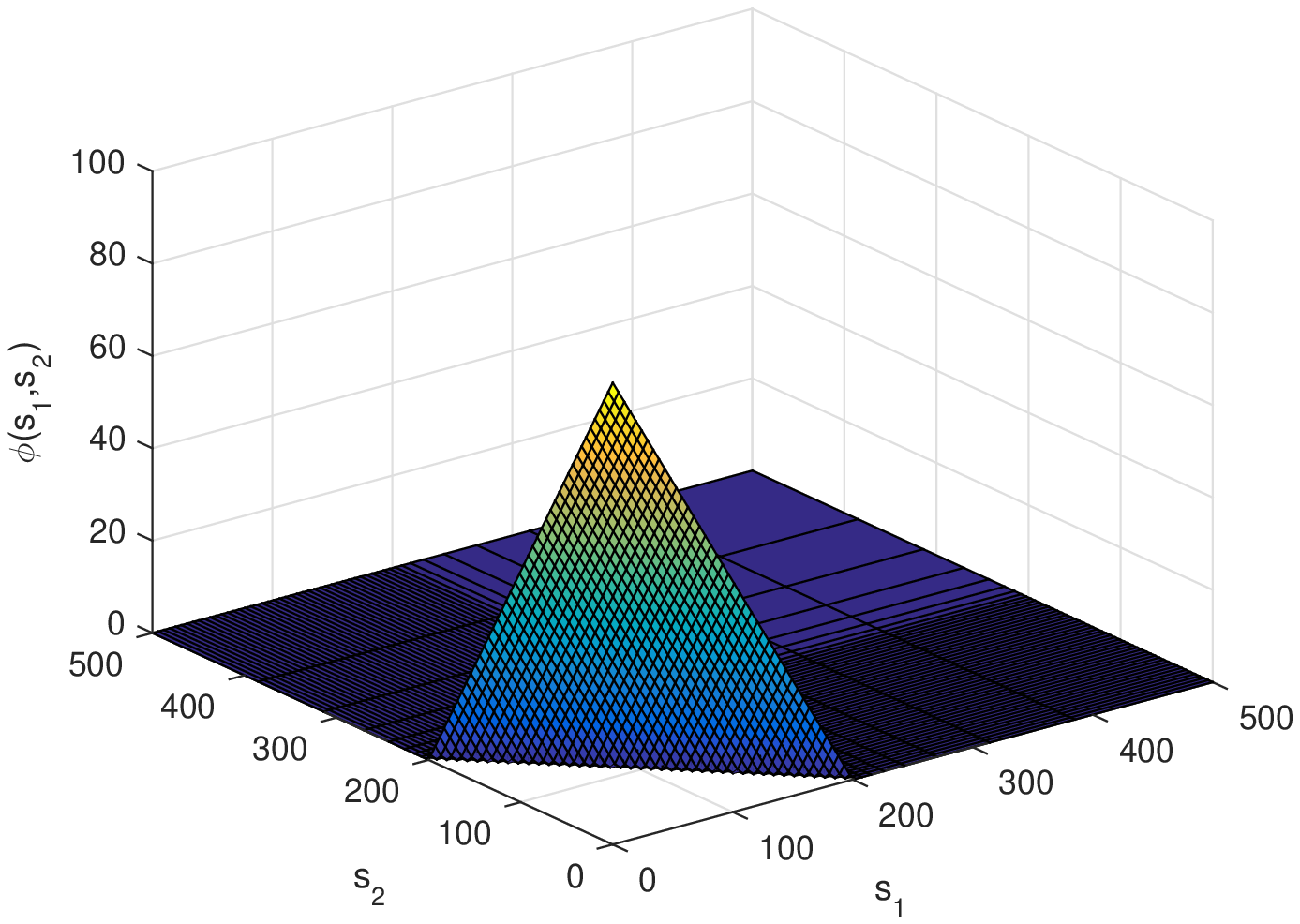}
    \caption{Payoff of the European put-on-the-min (left) and the European put-on-the-average (right) option if $K = 100$, $S_{\rm max} = 5K$.}
    \label{figpayoff}
\end{figure}

Clearly, the two payoff functions $\phi$ are continuous, but not everywhere continuously differentiable. 
For each option a specific nonuniform spatial grid $\{(s_{1,i},s_{2,j}): 0\leq i\leq m_1, 0 \leq j \leq m_2\}$ is chosen,
where the nonsmoothness of $\phi$ is taken into account.
The number of grid points is taken to be the same in both directions, $m_1 = m_2 = m$.
For the put-on-the-min option, relatively many points $s_{1,i}$ and $s_{2,j}$ are placed around the locations $s_1=K$ and 
$s_2=K$, respectively.  
For the put-on-the-average option, the payoff is nonsmooth on the line segment given by $s_1+s_2 = 2K$. 
Here the grid in each direction is taken to be uniform on the interval $[0,2K]$ and nonuniform outside, with relatively 
many points inside this interval. 
For both options, the nonuniform grids are generated by applying a smooth transformation to an artificial uniform grid, 
see e.g.~\cite{H17book}.
For Sets 1, 2, 3 we heuristically choose $S_{\rm max} = 5K, 30K, 50K$, respectively, in the case of the put-on-the-min
and $S_{\rm max} = 5K, 15K, 25K$, respectively, in the case of the put-on-the-average option.
As an illustration, Figure~\ref{figgrid} displays the two spatial grids if $m=50$, $K = 100$, $S_{\rm max} = 5K$.
In Figure~\ref{figpayoff} the graphs of the two payoff functions are shown.

The initial vector $V(0)=V^0$ is defined by pointwise evaluation of the payoff function $\phi$ at the spatial grid points, 
except at those points that lie closest to the set of nonsmoothness of $\phi$, where cell averaging is applied. 
For $k=1,2$ let
\begin{eqnarray*}
&s_{k,l+1/2} = \tfrac{1}{2}(s_{k,l}+s_{k,l+1})~~~~ &{\rm if}~~ 0\le l < m,\\
&h_{k,l+1/2} = s_{k,l+1/2}-s_{k,l-1/2} &{\rm if}~~ 0\le l \le m,
\end{eqnarray*}
with $s_{k,-1/2} = 0$ and $s_{k,m+1/2} = S_{\rm max}$.
Then we define \cite{H17book}
\[
V_{i,j}(0) = \frac{1}{h_{1,i+1/2}h_{2,j+1/2}}\int_{s_{1,i-1/2}}^{s_{1,i+1/2}}\int_{s_{2,j-1/2}}^{s_{2,j+1/2}} \phi(s_1,s_2) \dd s_2 \dd s_1
\]
whenever the cell 
\[
[s_{1,i-1/2}, s_{1,i+1/2})\times[s_{2,j-1/2}, s_{2,j+1/2})
\]
has a nonempty intersection with the set of points where $\phi$ is nonsmooth. 
For both payoffs under consideration, the above integral is readily calculated.

In the following we investigate, for all seven operator splitting methods formulated in Section~\ref{SecTime} for the semidiscrete 
system \eqref{ODEs}, the temporal discretization error at $t=T$ on a natural region of financial interest (ROI),
\begin{align}\label{temperror}
\widehat{E}^{ROI}(m,N) = \max \left\{|V^{N'}_{i,j}-V_{i,j}(T)|: \Delta t = T/N'~\textrm{and}~\tfrac{1}{2}K < s_{1,i},s_{2,j}< \tfrac{3}{2}K \right\}.
\end{align}
Clearly, the temporal discretization error is measured in the maximum norm, which is the most relevant norm in finance.
Numerical experiments reveal that the evaluation of the two-dimensional integral part is computationally the dominant 
part in each time step, compare also~\cite{STS14}.
In view of this, for a fair comparison of the different time-stepping methods, we consider the three schemes \eqref{CNFI}, 
\eqref{IETR}, \eqref{MCS} 
with $N'=N$ time steps and the four schemes \eqref{CNFE}, \eqref{CNAB}, \eqref{MCS2}, \eqref{SC2A} with $N' = 2N$ time steps, as the latter 
four schemes require only one evaluation of the integral part per time step, whereas the former three schemes employ two such 
evaluations.
A reference value for $V(T)$ has been computed by applying the MCS2 scheme with $3000$ time steps.

\begin{figure}[h!]
    \centering
    \hspace{-0.5cm}\includegraphics[trim={0cm 0cm 1cm 0cm},clip,
    width =0.51\textwidth]{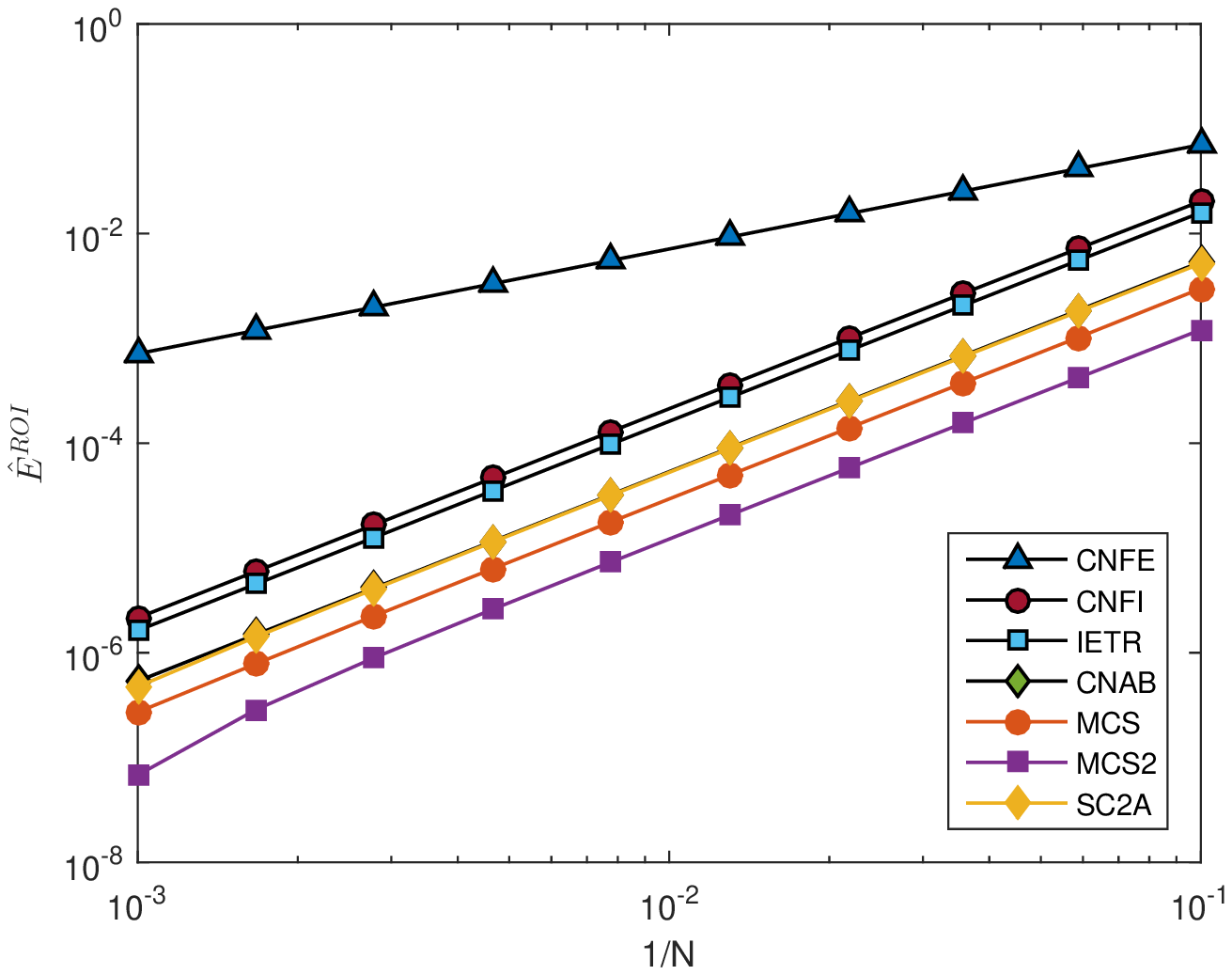}
    \includegraphics[trim={0cm 0cm 1cm 0cm},clip,
    width =0.51\textwidth]{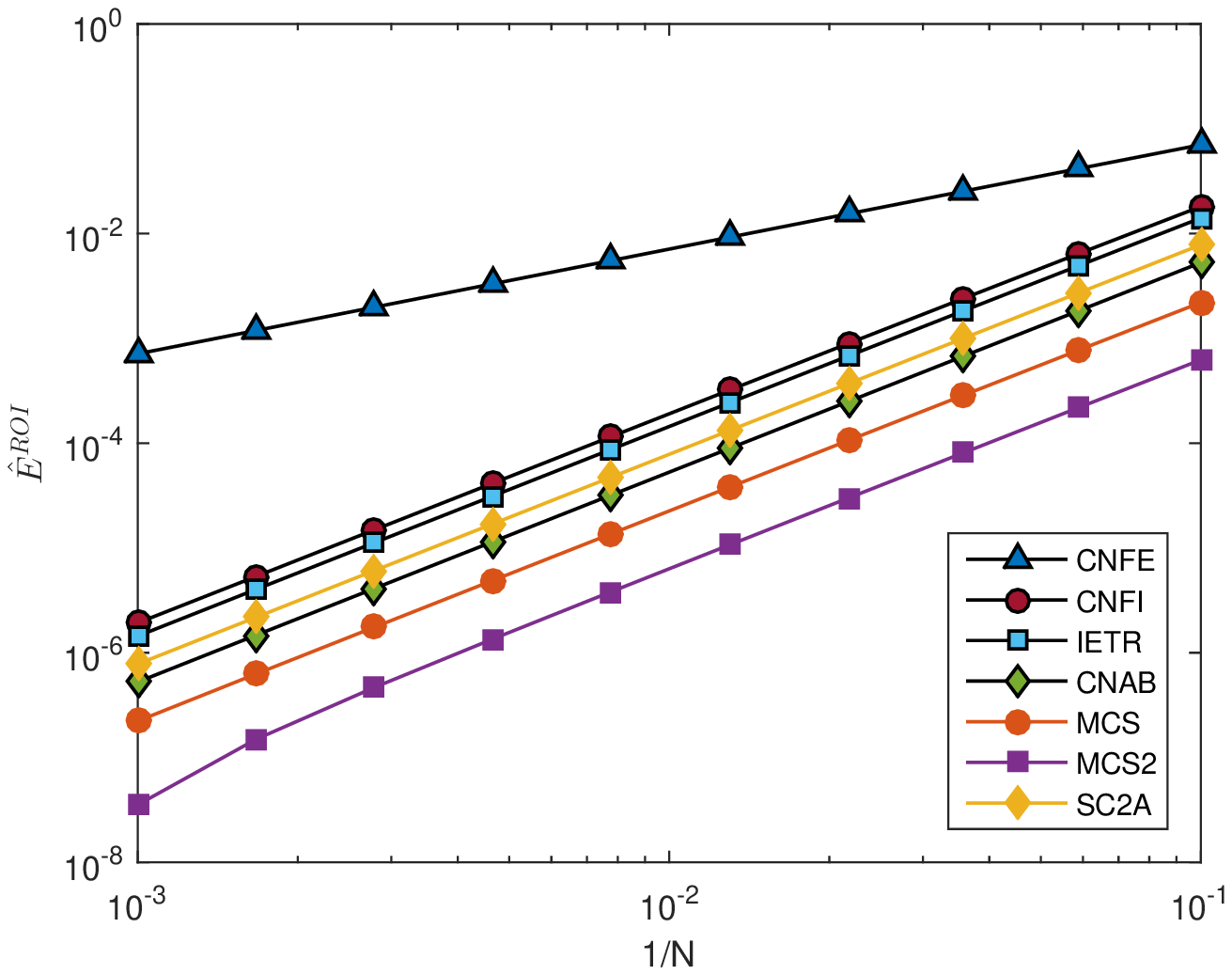}\\
    \hspace{-0.5cm}\includegraphics[trim={0cm 0cm 1cm 0cm},clip,
    width =0.51\textwidth]{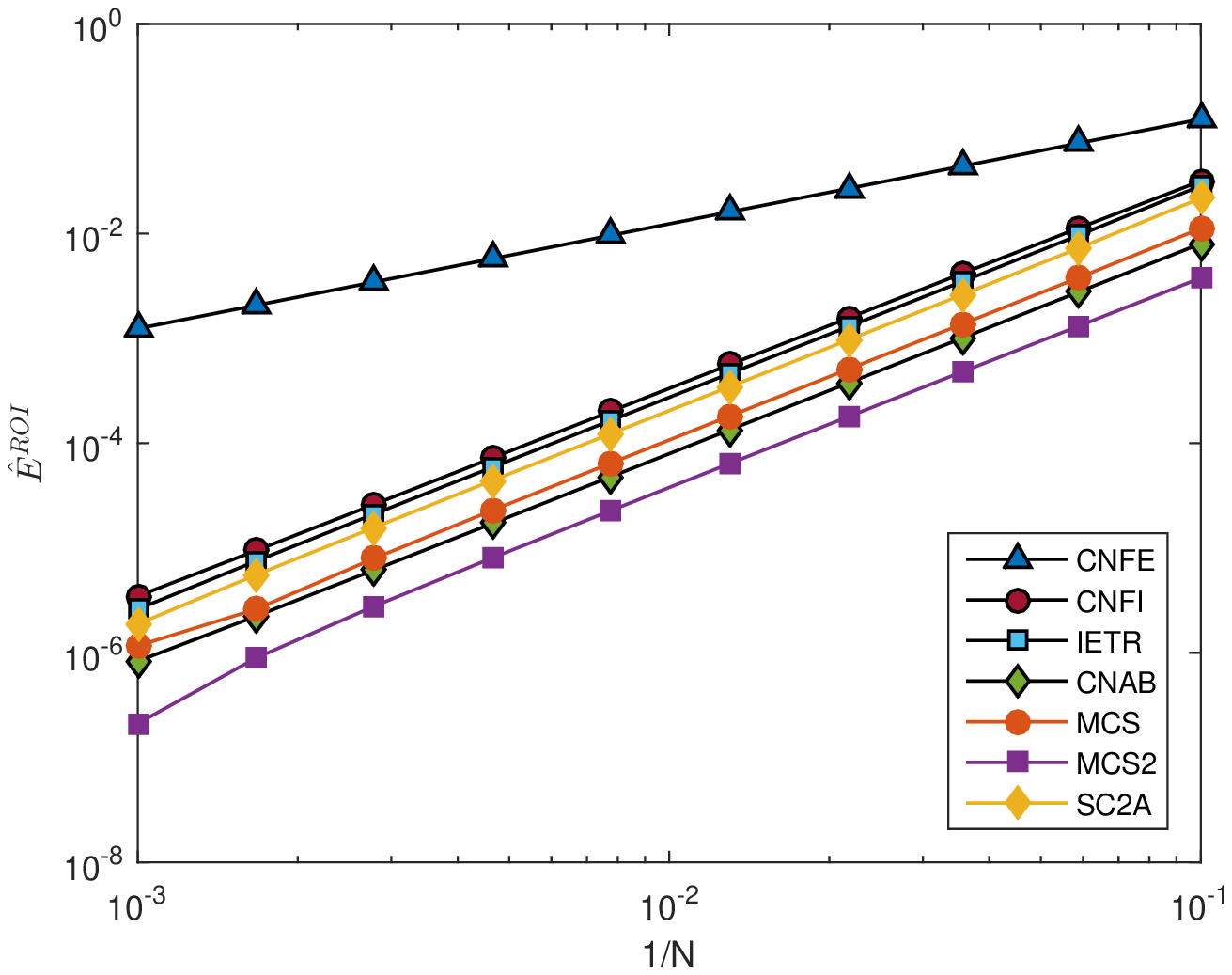}
    \includegraphics[trim={0cm 0cm 1cm 0cm},clip,
    width =0.51\textwidth]{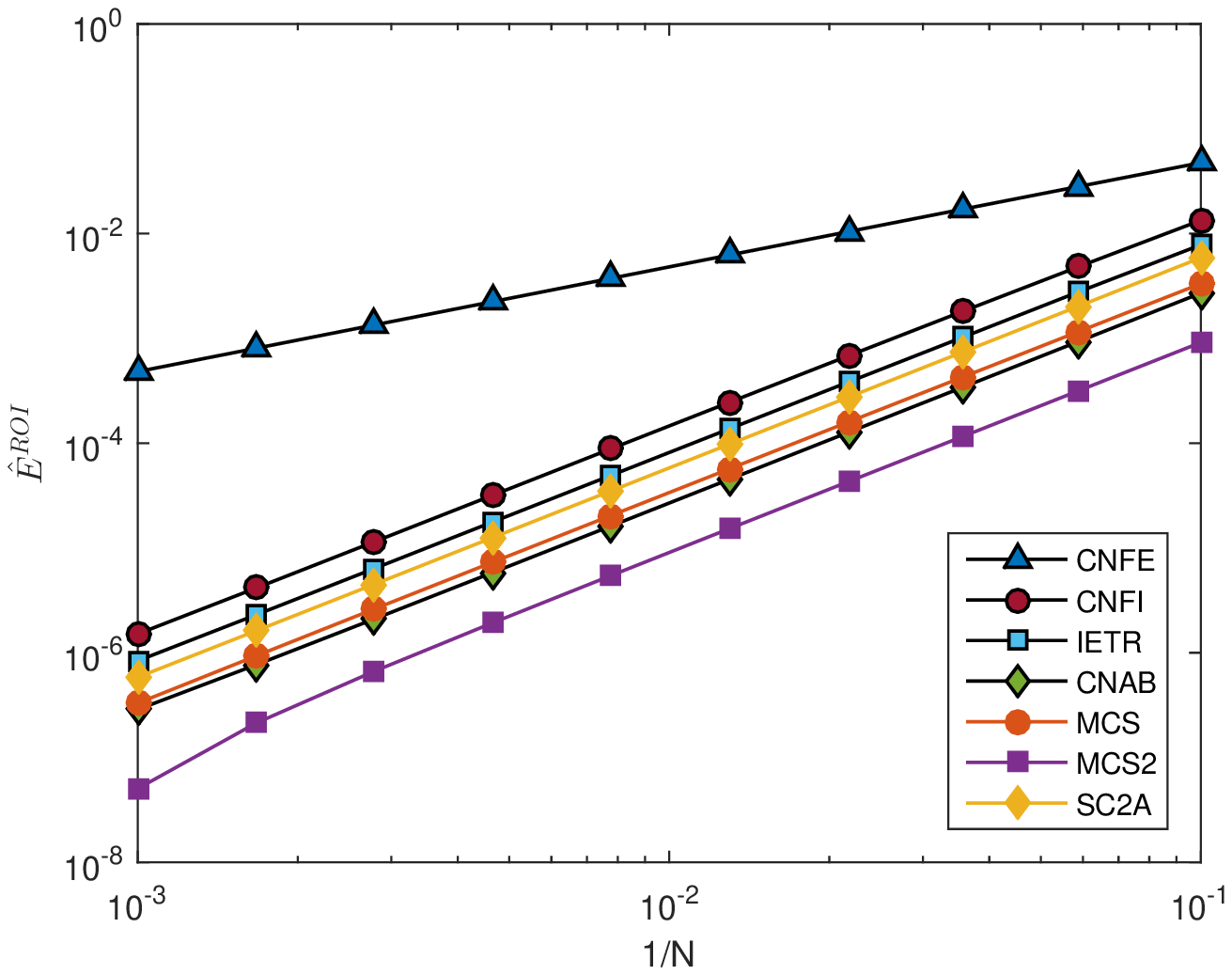}\\
    \hspace{-0.5cm}\includegraphics[trim={0cm 0cm 1cm 0cm},clip,
    width =0.51\textwidth]{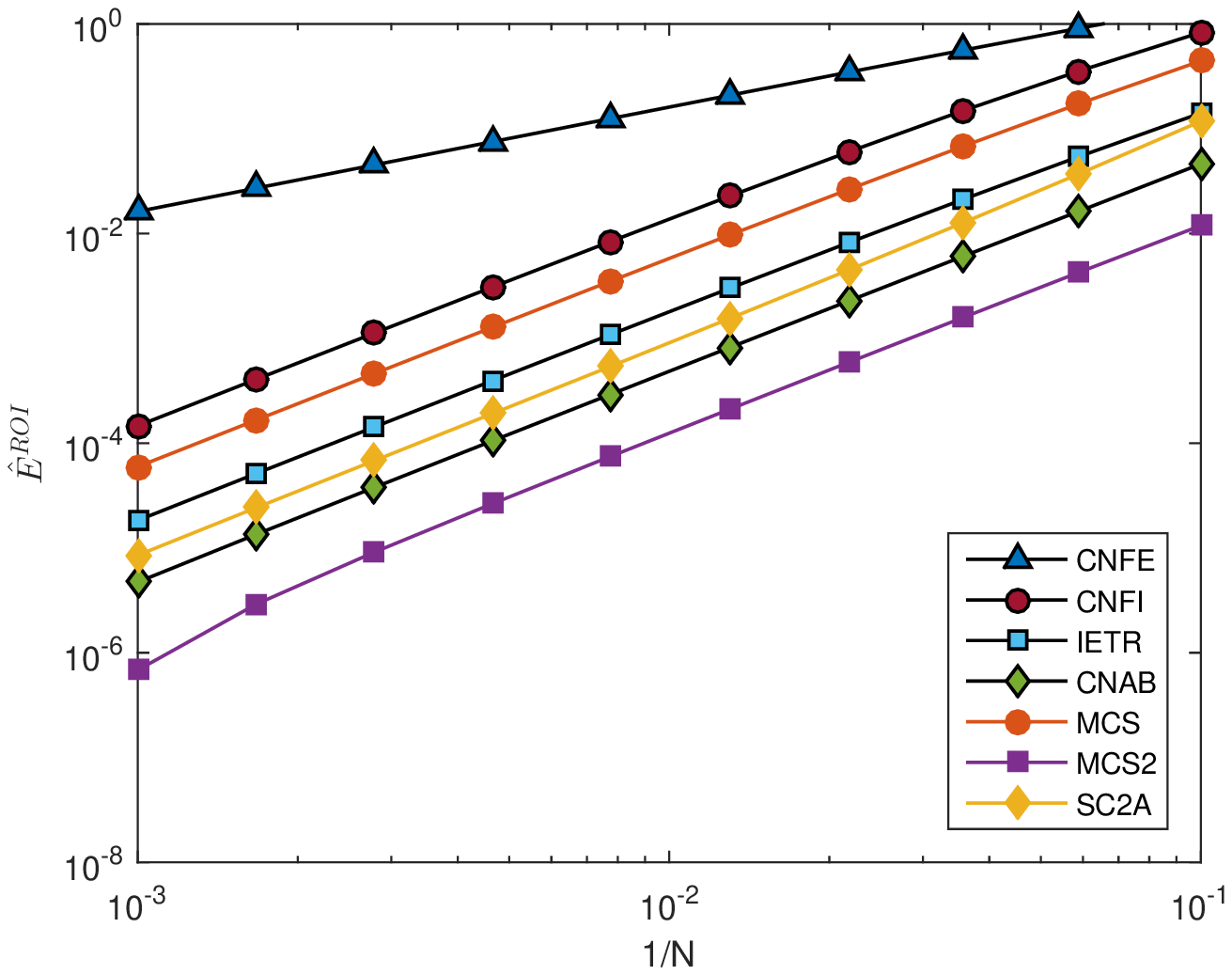}
    \includegraphics[trim={0cm 0cm 1cm 0cm},clip,
    width =0.51\textwidth]{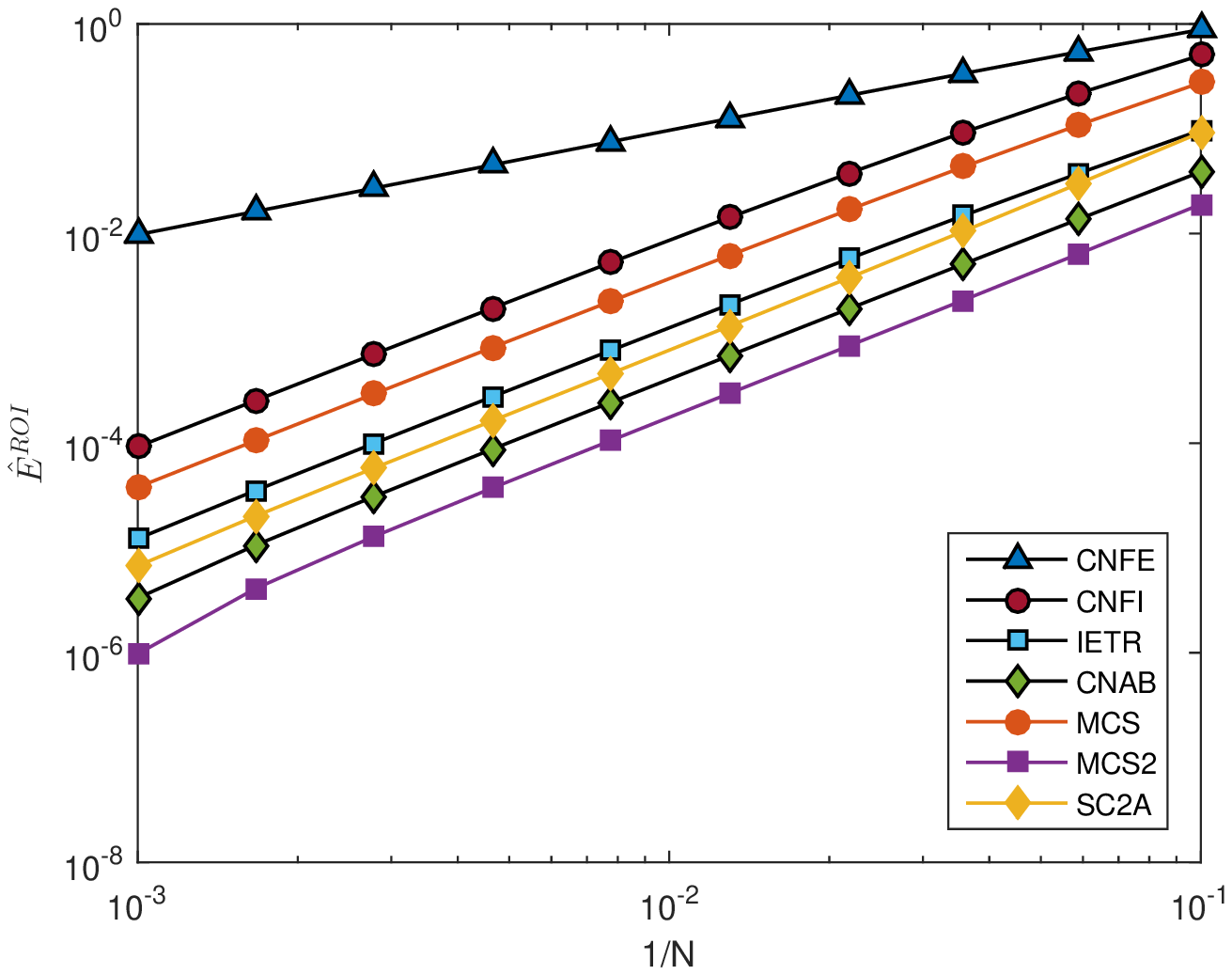}\\
    \caption{Temporal errors $\widehat{E}^{ROI}(150,N)$ of the seven operator splitting schemes in the case of the European put-on-the-min (left) 
    and the put-on-the-average (right) option under the two-asset Merton model with parameter Set~1 (top), Set~2 (mid) and Set~3 (bottom) from 
    Table~\ref{tabpars}. The CNFI, IETR, MCS schemes are applied with $\Delta t = T/N$ and the CNFE, CNAB, MCS2, SC2A schemes with $\Delta t = T/(2N)$.}
    \label{figConvtempset1}
\end{figure}

Figure~\ref{figConvtempset1} displays the temporal errors $\widehat{E}^{ROI}(m,N)$ for $m=150$ and a range of values $10 \le N \le 1000$ for
the put-on-the-min option (left column) and put-on-the-average option (right column) and the three parameter sets given in Table~\ref{tabpars}.
For all seven splitting methods, the positive result holds that the temporal errors always lie below a moderate value and show a regular, 
monotonically decreasing behaviour as $N$ increases.
The CNFE scheme \eqref{CNFE} shows an order of convergence just equal to one, as expected.
Each of the other six splitting schemes \eqref{CNFI}--\eqref{SC2A} attains a favourable order of convergence equal to two.
Comparing the different second-order methods, the MCS2 scheme yields the smallest temporal error constant in all experiments.
It is further interesting to notice that for Set~3, with more intensive jumps, each of the two-step schemes CNAB, MCS2, SC2A outperforms 
all of the one-step schemes under consideration.
The range of the temporal error constants for the different second-order methods is also the largest for this set.
We performed additional numerical experiments, choosing a coarser ($m=75$) as well as a finer ($m=250$) spatial grid. 
The obtained results were almost identical to those displayed in Figure~\ref{figConvtempset1}, indicating that the observed temporal 
convergence orders are valid in a favourable, stiff sense.
Recall that the CNFI scheme has been applied with a fixed number of $l=2$ iterations.
We also carried out experiments for this scheme with $l=3$ iterations, taking now $N' = \lceil 2N/3 \rceil$ time steps for 
a fair comparison. 
This did not lead to a significant improvement in accuracy, however, and did not alter the above conclusions.

\begin{figure}[h!]
    \centering
   \hspace{-0.5cm} \includegraphics[trim={0cm 0cm 1cm 0.5cm},clip,width =0.58\textwidth]{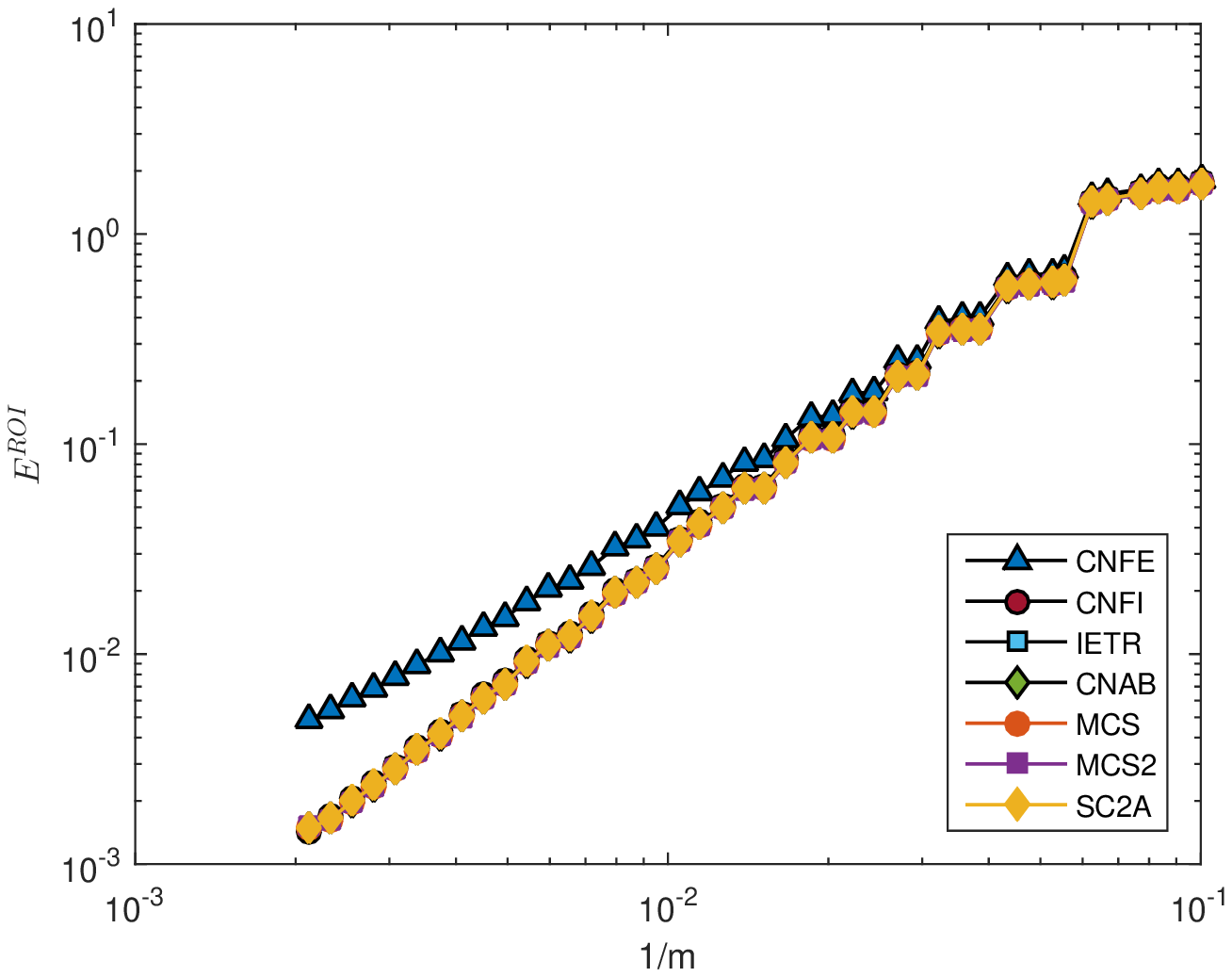}\\
   \hspace{-0.5cm} \includegraphics[trim={0cm 0cm 1cm 0.5cm},clip,width =0.58\textwidth]{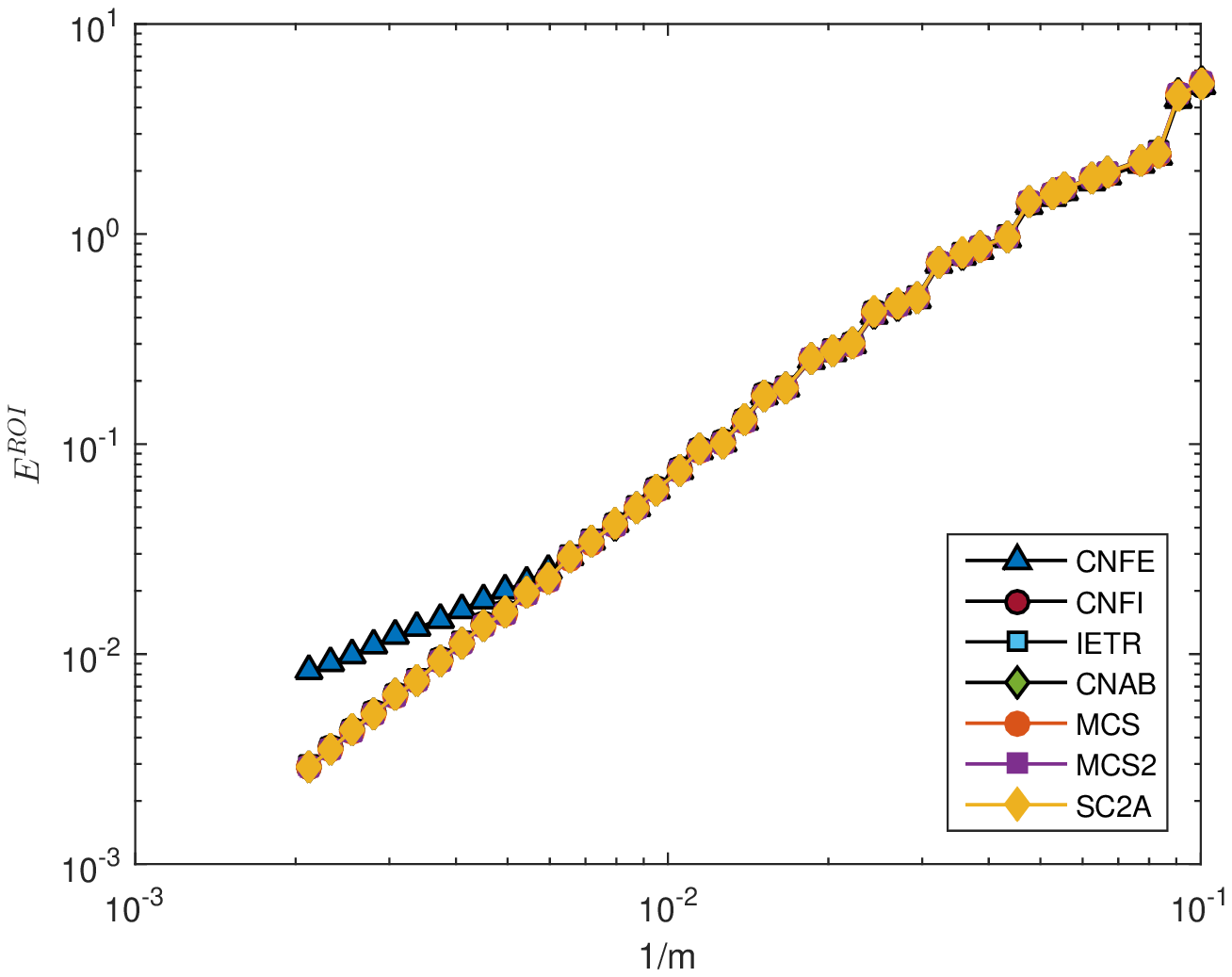}\\
    \hspace{-0.5cm}\includegraphics[trim={0cm 0cm 1cm 0.5cm},clip,width =0.58\textwidth]{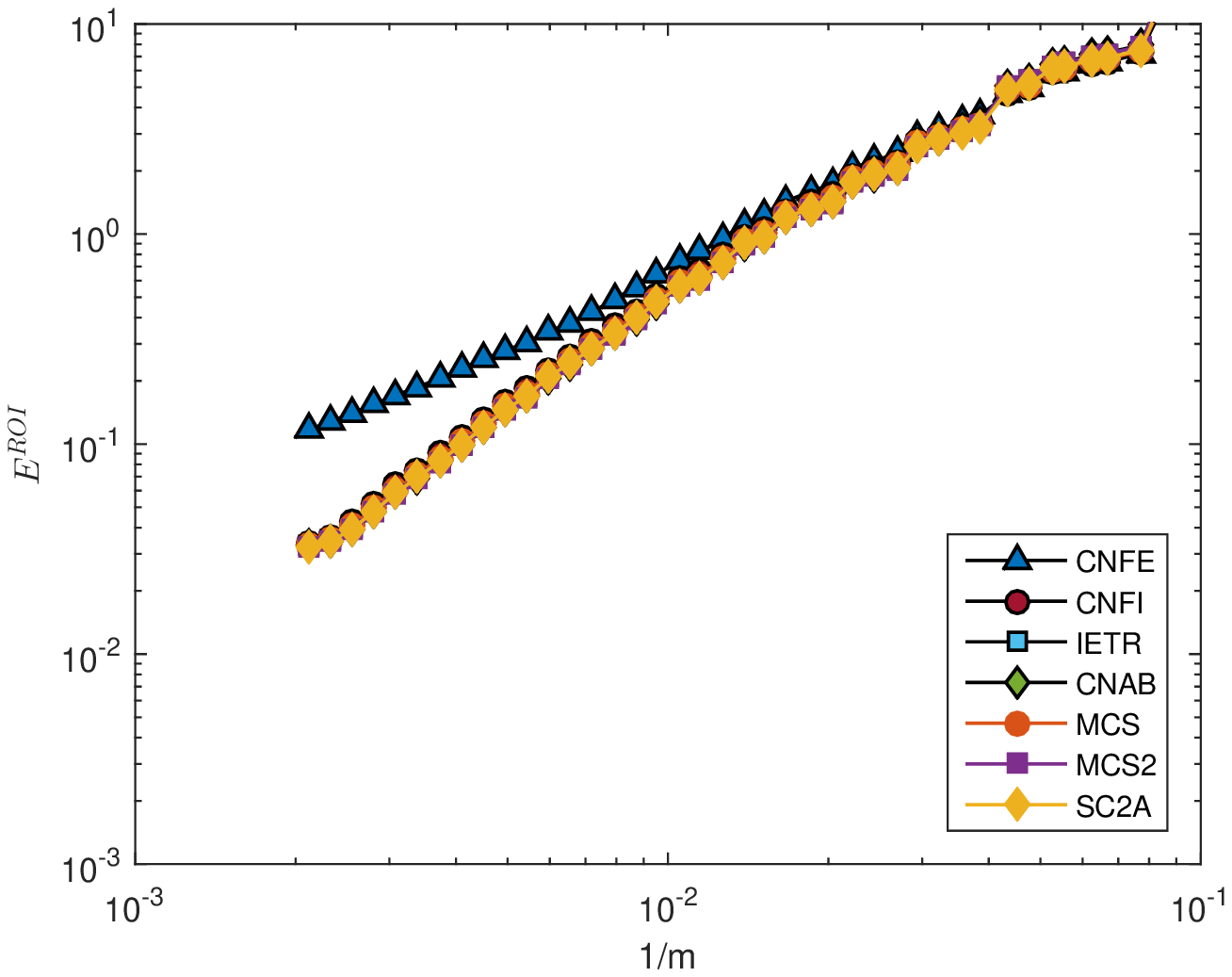}
    \caption{Total errors $E^{ROI}(m,N)$ in the case of the European put-on-the-min option under the two-asset Merton model with 
    parameter Set~1 (top), Set~2 (mid) and Set~3 (bottom) from Table~\ref{tabpars}. The CNFI, IETR, MCS schemes are applied with $\Delta t = T/N$ 
    and the CNFE, CNAB, MCS2, SC2A schemes with $\Delta t = T/(2N)$ and $N=\lceil m/3 \rceil$.}
    \label{figConvPutonMin}
\end{figure}

A semi-closed analytical formula has recently been derived in Boen \cite{B18} for the value of a European put-on-the-min option 
under the two-asset Merton model. This formula is stated in Appendix~\ref{SecExact}. 
Hence, in this case, we can also study the total discretization error at $t=T$ on the region of interest,
\begin{align} \label{globerror}
    E^{ROI}(m,N) = \max \left\{|V^{N'}_{i,j}-v(s_{1,i},s_{2,j},T)|: \Delta t = T/N'~\textrm{and}~\tfrac{1}{2}K < s_{1,i},s_{2,j}< \tfrac{3}{2}K \right\}.
\end{align}
Figure~\ref{figConvPutonMin} displays the total errors $E^{ROI}(m,N)$ for a range of 
values $10\le m \le 500$ and the three parameter sets given in Table \ref{tabpars}, where we have chosen $N=\lceil m/3 \rceil$ and 
$N'$ as above.
Clearly, with all splitting methods except CNFE, a second-order convergence behaviour is found, which is as desired.
The total errors obtained with these six time-stepping methods are visually almost identical, indicating that the spatial 
discretization error dominates in these experiments.


\section{Conclusion}\label{SecConc}
In this paper we have studied seven different operator splitting schemes when applied to the two-dimensional Merton PIDE with a 
focus on IMEX and ADI methods. 
Each of the considered schemes conveniently treats the nonstiff, nonlocal integral part in an explicit fashion. 
Any matrix-vector product arising from the semidiscretization of the integral part is computed efficiently by means of FFT using 
Algorithm~\ref{algo1}. 
Through ample numerical experiments, we have investigated the convergence behaviour of the seven schemes and examined their relative 
performance. 
Except for the first-order CNFE scheme, all of the considered splitting schemes showed a stiff order of temporal convergence equal 
to two, which is as desired. 
The two-step MCS2 scheme with $\theta = \frac{1}{3}$ stood out as yielding the smallest temporal error constant in all the conducted 
experiments. 
Together with the positive results recently derived in \cite{HT18}, we therefore recommend this splitting scheme for the 
efficient and stable temporal discretization of two-dimensional PIDEs in financial option valuation.

\section*{Acknowledgements}
The authors acknowledge the support of the Research Fund (BOF) of the University of Antwerp (41/FA070300/3/FFB150337).


\clearpage
\section*{Figures}

\clearpage
\appendix
\section{The block-Toeplitz structure of $\overline{A}^{(J)}$}\label{appToep}
As an illustration, if $M_1 = M_2 =4$, then $\overline{F}$ in the matrix $\overline{A}^{(J)} = \lambda \of$ is of the following form:
\begin{figure}[h!]
\centering
\begin{tikzpicture}[ ]
\useasboundingbox (-6.5,-4.6) rectangle (6.5,5);
    \scope[transform canvas={scale=.9}]
         \matrix (m)[matrix of math nodes,left delimiter=(,right delimiter=)]
	{
	\of_{0,0} & \of_{1,0} & \of_{2,0} & \of_{3,0} & \of_{0,1} & \of_{1,1} & \of_{2,1} & \of_{3,1} & \of_{0,2} & \of_{1,2} & \of_{2,2} & \of_{3,2} & \of_{0,3} & \of_{1,3} & \of_{2,3} & \of_{3,3}\\
	\of_{\minn1,0} & \of_{0,0} & \of_{1,0} & \of_{2,0} & \of_{\minn1,1} & \of_{0,1} & \of_{1,1} & \of_{2,1} & \of_{\minn1,2} & \of_{0,2} & \of_{1,2} & \of_{2,2} & \of_{\minn1,3} & \of_{0,3} & \of_{1,3} & \of_{2,3} \\
	\of_{\minn2,0} & \of_{\minn1,0} & \of_{0,0} & \of_{1,0} & \of_{\minn2,1} & \of_{\minn1,1} & \of_{0,1} & \of_{1,1} & \of_{\minn2,2} & \of_{\minn1,2} & \of_{0,2} & \of_{1,2} & \of_{\minn2,3} & \of_{\minn1,3} & \of_{0,3} & \of_{1,3}\\
	\of_{\minn3,0} & \of_{\minn2,0} & \of_{\minn1,0} & \of_{0,0} & \of_{\minn3,1} & \of_{\minn2,1} & \of_{\minn1,1} & \of_{0,1} & \of_{\minn3,2} & \of_{\minn2,2} & \of_{\minn1,2} & \of_{0,2} & \of_{\minn3,3} & \of_{\minn2,3} & \of_{\minn1,3} & \of_{0,3} \\
	\of_{0,\minn1} & \of_{1,\minn1} & \of_{2,\minn1} & \of_{3,\minn1} & \of_{0,0} & \of_{1,0} & \of_{2,0} & \of_{3,0} & \of_{0,1} & \of_{1,1} & \of_{2,1} & \of_{3,1} & \of_{0,2} & \of_{1,2} & \of_{2,2} & \of_{3,2} \\
	\of_{\minn1,\minn1} & \of_{0,\minn1} & \of_{1,\minn1} & \of_{2,\minn1} &\of_{\minn1,0} & \of_{0,0} & \of_{1,0} & \of_{2,0} & \of_{\minn1,1} & \of_{0,1} & \of_{1,1} & \of_{2,1} & \of_{\minn1,2} & \of_{0,2} & \of_{1,2} & \of_{2,2}  \\
	\of_{\minn2,\minn1} & \of_{\minn1,\minn1} & \of_{0,\minn1} & \of_{1,\minn1} &\of_{\minn2,0} & \of_{\minn1,0} & \of_{0,0} & \of_{1,0} & \of_{\minn2,1} & \of_{\minn1,1} & \of_{0,1} & \of_{1,1} & \of_{\minn2,2} & \of_{\minn1,2} & \of_{0,2} & \of_{1,2}\\
	\of_{\minn3,\minn1} &\of_{\minn2,\minn1} &  \of_{\minn1,\minn1} & \of_{0,\minn1} &\of_{\minn3,0} & \of_{\minn2,0} & \of_{\minn1,0} & \of_{0,0} & \of_{\minn3,1} & \of_{\minn2,1} & \of_{\minn1,1} & \of_{0,1} & \of_{\minn3,2} & \of_{\minn2,2} & \of_{\minn1,2} & \of_{0,2}\\
	\of_{0,\minn2} & \of_{1,\minn2} & \of_{2,\minn2} & \of_{3,\minn2} & \of_{0,\minn1} & \of_{1,\minn1} & \of_{2,\minn1} & \of_{3,\minn1} & \of_{0,0} & \of_{1,0} & \of_{2,0} & \of_{3,0} & \of_{0,1} & \of_{1,1} & \of_{2,1} & \of_{3,1}\\
	\of_{\minn1,\minn2} &\of_{0,\minn2} & \of_{1,\minn2} & \of_{2,\minn2} & \of_{\minn1,\minn1} & \of_{0,\minn1} & \of_{1,\minn1} & \of_{2,\minn1} &\of_{\minn1,0} & \of_{0,0} & \of_{1,0} & \of_{2,0} & \of_{\minn1,1} & \of_{0,1} & \of_{1,1} & \of_{2,1}  \\
	\of_{\minn2,\minn2} & \of_{\minn1,\minn2} &\of_{0,\minn2} & \of_{1,\minn2} &\of_{\minn2,\minn1} & \of_{\minn1,\minn1} & \of_{0,\minn1} & \of_{1,\minn1} &\of_{\minn2,0} & \of_{\minn1,0} & \of_{0,0} & \of_{1,0} & \of_{\minn2,1} & \of_{\minn1,1} & \of_{0,1} & \of_{1,1}\\
	\of_{\minn3,\minn2} &\of_{\minn2,\minn2} & \of_{\minn1,\minn2} &\of_{0,\minn2} &\of_{\minn3,\minn1} &\of_{\minn2,\minn1} &  \of_{\minn1,\minn1} & \of_{0,\minn1} &\of_{\minn3,0} & \of_{\minn2,0} & \of_{\minn1,0} & \of_{0,0} & \of_{\minn3,1} & \of_{\minn2,1} & \of_{\minn1,1} & \of_{0,1} \\
	\of_{0,\minn3} & \of_{1,\minn3} & \of_{2,\minn3} & \of_{3,\minn3} & \of_{0,\minn2} & \of_{1,\minn2} & \of_{2,\minn2} & \of_{3,\minn2} & \of_{0,\minn1} & \of_{1,\minn1} & \of_{2,\minn1} & \of_{3,\minn1} & \of_{0,0} & \of_{1,0} & \of_{2,0} & \of_{3,0} \\
	\of_{\minn1,\minn3} & \of_{0,\minn3} & \of_{1,\minn3} & \of_{2,\minn3} &\of_{\minn1,\minn2} &\of_{0,\minn2} & \of_{1,\minn2} & \of_{2,\minn2} & \of_{\minn1,\minn1} & \of_{0,\minn1} & \of_{1,\minn1} & \of_{2,\minn1} &\of_{\minn1,0} & \of_{0,0} & \of_{1,0} & \of_{2,0}  \\
	\of_{\minn2,\minn3} & \of_{\minn1,\minn3} & \of_{0,\minn3} & \of_{1,\minn3} & \of_{\minn2,\minn2} & \of_{\minn1,\minn2} &\of_{0,\minn2} & \of_{1,\minn2} &\of_{\minn2,\minn1} & \of_{\minn1,\minn1} & \of_{0,\minn1} & \of_{1,\minn1} &\of_{\minn2,0} & \of_{\minn1,0} & \of_{0,0} & \of_{1,0} \\
	\of_{\minn3,\minn3} & \of_{\minn2,\minn3} & \of_{\minn1,\minn3} & \of_{0,\minn3} & \of_{\minn3,\minn2} &\of_{\minn2,\minn2} & \of_{\minn1,\minn2} &\of_{0,\minn2} &\of_{\minn3,\minn1} &\of_{\minn2,\minn1} &  \of_{\minn1,\minn1} & \of_{0,\minn1} &\of_{\minn3,0} & \of_{\minn2,0} & \of_{\minn1,0} & \of_{0,0} \\
	};

\draw[scale=0.9](m-1-4.north east)--(m-16-4.south east);
\draw[scale=0.9](m-1-8.north east)--(m-16-8.south east);
\draw[scale=0.9](m-1-12.north east)--(m-16-12.south east);

\draw[scale=0.9](m-4-1.south west)--(m-4-16.south east);
\draw[scale=0.9](m-8-1.south west)--(m-8-16.south east);
\draw[scale=0.9](m-12-1.south west)--(m-12-16.south east);
  
 \endscope

\end{tikzpicture}
\end{figure}

\hspace{0.3cm}\\
This is a block-Toeplitz matrix with Toeplitz blocks $\of^1_d$, ${d=-3,\ldots,3}$:
\[
\overline{F} = \begin{pmatrix}
\of^1_0 & \of^1_1 & \of^1_2 & \of^1_3 \\
\of^1_{\minn1} & \of^1_0 & \of^1_1 & \of^1_2 \\
\of^1_{\minn2} & \of^1_{\minn1} & \of^1_0 & \of^1_1 \\
\of^1_{\minn3} & \of^1_{\minn2} & \of^1_{\minn1} & \of^1_0
\end{pmatrix}.
\]

\clearpage
\section{Semi-closed analytical formula for put-on-the-min option}\label{SecExact}
Let $M(x_1,x_2,\Sigma)$ denote the bivariate normal cumulative distribution function, evaluated at $(x_1,x_2)$, with mean $(0,0)^\top$ and covariance matrix $\Sigma$.
Under the two-asset Merton jump-diffusion model, the value of a European put-on-the-min option at inception is given in Boen \cite{B18}:
\begin{align*}\label{exactprice}
    v(S^{(1)}_0,S^{(2)}_0,T) = \sum_{n=0}^{\infty} &e^{-\lambda T}\frac{(\lambda T)^n}{n!} \left[e^{-rT}K-e^{-rT}KM(b_1,b_2,\Sigma_3)-\right.\\
    &\left.S^{(1)}_0e^{-\lambda\kappa_1T + n\gamma_1+\half n \delta_1^2}M(d_{11},d_1,\Sigma_1) - S^{(2)}_0e^{-\lambda\kappa_2T + n\gamma_2 + \half n \delta_2^2}M(d_{22},d_2,\Sigma_2)\right],\nonumber
\end{align*}
where
\[
b_i = \frac{\ln\left(\frac{S^{(i)}_0}{K}\right)+(r-\half\sigma_i^2-\lambda\kappa_i)T+n\gamma_i}{\sqrt{T}\sqrt{\sigma_i^2+\frac{n\delta_i^2}{T}}}, \quad d_i = -b_i - \sqrt{T}\sqrt{\sigma_i^2+\frac{n\delta_i^2}{T}}\quad (i=1,2),
\]
\[
d_{11} = \frac{\ln\left(\frac{S^{(2)}_0}{S^{(1)}_0}\right) +\left(-\half\sigma^2+\lambda(\kappa_1-\kappa_2)\right)T-n(\gamma_1-\gamma_2+\delta_1^2-\rhoh\delta_1\delta_2)}{\sqrt{T}\sqrt{\sigma^2+\frac{n\delta^2}{T}}},\quad d_{22} = -d_{11}-\sqrt{T}\sqrt{\sigma^2+\frac{n\delta^2}{T}},
\]
\[ 
\Sigma_i= \begin{pmatrix}
1 & \rho_i\\
\rho_i & 1
\end{pmatrix}\quad (i=1,2,3),
\]
with
\[
\sigma^2 = \sigma_1^2-2\rho\sigma_1\sigma_2+\sigma_2^2, \quad \delta^2 = \delta_1^2-2\rhoh\delta_1\delta_2+\delta_2^2,
\]
\[
\rho_i = \frac{\sqrt{\sigma_i^2+\frac{n\delta_i^2}{T}}}{\sqrt{\sigma^2+\frac{n\delta^2}{T}}}-\frac{\rho\sigma_1\sigma_2+\frac{n\rhoh\delta_1\delta_2}{T}}{\sqrt{\left(\sigma^2+\frac{n\delta^2}{T}\right)\left(\sigma_i^2+\frac{n\delta_i^2}{T}\right)}}\quad (i=1,2),
\]
and
\[
\rho_3 = \frac{\rho\sigma_1\sigma_2+\frac{n\rhoh\delta_1\delta_2}{T}}{\sqrt{\left(\sigma_1^2+\frac{n\delta_1^2}{T}\right)\left(\sigma_2^2+\frac{n\delta_2^2}{T}\right)}}.
\]
\clearpage


\bibliographystyle{plain}
\bibliography{bib_2DMerton_splitting}

\end{document}